\documentclass[10pt]{amsart}
\usepackage[cp1251]{inputenc}
\usepackage[german,english]{babel}
\usepackage{amsmath}
\usepackage{amssymb}
\usepackage{amsfonts}

\begin{document}
\renewcommand{\refname}{References}

\thispagestyle{empty}

\title[Application of Multiple Fourier--Legendre Series]
{Application of Multiple Fourier--Legendre Series to  
Implementation of Strong Exponential Milstein and 
Wagner--Platen Methods for Non-Commutative Semilinear Stochastic 
Partial Differential Equations}
\author[D.F. Kuznetsov]{Dmitriy F. Kuznetsov}
\address{Dmitriy Feliksovich Kuznetsov
\newline\hphantom{iii} Peter the Great Saint-Petersburg Polytechnic University,
\newline\hphantom{iii} Institute of Applied Mathematics and Mechanics,
\newline\hphantom{iii} Department of Mathematics,
\newline\hphantom{iii} Polytechnicheskaya ul., 29,
\newline\hphantom{iii} 195251, Saint-Petersburg, Russia}
\email{sde\_kuznetsov@inbox.ru}
\thanks{\sc Mathematics Subject Classification: 60H05, 60H10, 42B05, 42C10}
\thanks{\sc Keywords: Non-commutative semilinear
Stochastic partial differential equation,
Multiplicative trace class noise,
Infinite-dimensional $Q$-Wiener process, 
Iterated It\^{o} stochastic integral,
Generalized multiple Fourier series,
Multiple Fourier--Legendre series, Exponential Milstein scheme,
Exponential Wagner--Platen scheme,
Legendre polynomials, Mean-square approximation, Expansion.}

\maketitle {\small
\begin{quote}
\noindent{\sc Abstract.} 
The article is devoted to 
the application of multiple Fourier--Legendre series to  
implementation of strong exponential Milstein and 
Wagner--Platen methods for non-commutative semilinear stochastic 
partial differential equations with multiplicative 
trace class noise. These methods have 
strong orders of convergence $1.0-\varepsilon$ and $1.5-\varepsilon$
correspondingly (here $\varepsilon$ is an arbitrary small
positive real number) with respect to the temporal discretization.
The theorem on mean-square convergence
of approximations of iterated stochastic integrals
of multiplicities 1 to 3 with respect to the infinite-dimensional
$Q$-Wiener process is formulated and proved.
The results of the article can be applied to
implementation 
of exponential Milstein and 
Wagner--Platen methods for non-commutative semilinear stochastic 
partial differential equations
with multiplicative 
trace class noise.
\medskip
\end{quote}
}

\vspace{5mm}

\section{Introduction}

\vspace{5mm}

It is well-known that
one of the effective approaches to the construction of high-order
strong numerical methods (with respect to the temporal 
discretization) for stochastic partial differential
equations (SPDEs) is based
on the Taylor formula in Banach spaces and exponential formula for 
the mild solution of 
SPDEs \cite{1}-\cite{6}. 
A significant step in this direction 
was made in \cite{2}, \cite{3}, where the exponential 
Milstein and 
Wagner--Platen methods for semilinear SPDEs 
were constructed.
Under the appropriate conditions \cite{2}, \cite{3} these methods
have strong orders of convergence $1.0-\varepsilon$ and $1.5-\varepsilon$
correspondingly (where $\varepsilon$
is an arbitrary small posilive real number)
with respect to the temporal discretization.

An important feature of the mentioned numerical methods is the presence 
in them of the so-called iterated stochastic integrals with respect 
to the infinite-dimensional $Q$-Wiener process \cite{7}.
The problem of numerical modeling of these stochastic integrals
was solved in \cite{2}, \cite{3} for the case when special 
commutativity conditions are fulfilled.

If the mentioned commutativity conditions are not fulfilled,
which often corresponds to SPDEs in numerous applications, the 
numerical simulation of 
iterated stochastic integrals with respect 
to the infinite-dimensional $Q$-Wiener process becomes 
much more difficult. 
Note that the exponential Milstein scheme \cite{2} contains the
iterated stochastic integrals of multiplicities 1 and 2 with respect 
to the infinite-dimensional $Q$-Wiener process
and the exponential Wagner--Platen scheme \cite{3}
contains the mentioned stochastic integrals of multiplicities
1 to 3.
In \cite{8} two methods of the mean-square approximation of iterated
stochastic integrals from the Milstein scheme \cite{2}
have been considered. Note that the mean-square error of approximation
of these stochastic integrals
consists of two components \cite{8}. 
The first component is related with the finite-dimentional
approximation of the infinite-dimentional $Q$-Wiener process
while the second one is connected with the approximation
of iterated It\^{o} stochastic integrals with respect to the
scalar standard Brownian motions.
In the author's publications \cite{9}, \cite{arxiv-20}, \cite{24a}, \cite{24aa} 
the problem of the mean-square approximation 
of iterated stochastic integrals 
with respect to the infinite-dimensional $Q$-Wiener process
in the sense of second component of approximation error (see above)
has been solved 
for an arbitrary multiplicity $k$ ($k\in\mathbb{N}$)
of stochastic integrals.
More precisely, in \cite{9}, \cite{arxiv-20} (also see \cite{24a}, \cite{24aa}) the method of 
generalized multiple Fourier series \cite{11}-\cite{arxiv-25} for the
approximation of iterated It\^{o} stochastic integrals with respect to the 
scalar standard Brownian motions was adapted for  
iterated stochastic integrals with respect to the 
infinite-dimensional $Q$-Wiener process 
(in the sense of second component of the approximation error).

In this article, we extend the method for estimating the first 
component of approximation error
from \cite{8} for iterated stochastic integrals 
of multiplicities 1 to 3 with respect to the 
infinite-dimensional $Q$-Wiener process. 
In addition, we combine the obtained 
results with the results from \cite{9}, \cite{arxiv-20}. Thus, the results 
of the article can be applied to the implementation of exponential Milstein
and Wagner--Platen schemes for semilinear SPDEs 
with multiplicative trace class noise and
without the conditions of commutativity for SPDEs.

\vspace{5mm}

\section{Exponential Milstein and Wagner--Platen Numerical Schemes
for Non-Commutative Semilinear SPDEs}

\vspace{5mm}

Let $U, H$ be separable $\mathbb{R}$-Hilbert spaces and
$L_{HS}(U,H)$ be a space of Hilbert--Schmidt operators mapping from $U$ to $H.$
Let $(\Omega, {\bf F}, \sf{P})$ be a probability space 
with a normal filtration $\{{\bf F}_t, t\in [0, \bar{T}]\}$
\cite{7}, let ${\bf W}_t$ be an $U$-valued $Q$-Wiener process 
with respect to $\{{\bf F}_t, t\in [0, \bar{T}]\},$  
which has a covariance trace class operator $Q\in L(U)$. 
Here $L(U)$ denotes all bounded linear operators
mapping from $U$ to $U.$ Consider an $\mathbb{R}$-Hilbert space 
$U_0=Q^{1/2}(U)$ with a scalar product 

$$
\left\langle u, w\right\rangle_{U_0}=\left\langle 
Q^{-1/2}u, Q^{-1/2}w\right\rangle_{U}
$$

\vspace{3mm}
\noindent
for all $u, w\in U_0.$

Consider the semilinear parabolic 
SPDE 

\vspace{1mm}
\begin{equation}
\label{xx1}
dX_t = \left(A X_t + F(X_t)\right)dt + B(X_t)d{\bf W}_t,\ \ \
X_0=\xi,\ \ \ t\in [0, \bar{T}],
\end{equation}

\vspace{5mm}
\noindent
where nonlinear operators $F,$ $B$ ($F:$ $H\rightarrow H$, 
$B:$ $H\rightarrow L_{HS}(U_0,H)$), linear operator
$A:$ $D(A)\subset H\rightarrow H$
as well as the initial value $\xi$ 
are assumed to satisfy the conditions
of existence and uniqueness of
the SPDE (\ref{xx1}) mild solution (see \cite{3}, Assumptions A1--A4).

It is well-known \cite{6} that Assumptions A1--A4 \cite{3}
guarantee the existence and uniqueness (up to modifications) of 
the mild solution 
$X_t: [0, \bar T]\times \Omega \rightarrow H $
of the SPDE (\ref{xx1}) 

\vspace{1mm}
\begin{equation}
\label{mild}
X_t={\rm exp}(At)\xi +\int\limits_0^t {\rm exp}(A(t-\tau))F(X_\tau)d\tau+
\int\limits_0^t {\rm exp}(A(t-\tau))B(X_\tau)d{\bf W}_{\tau}
\end{equation}

\vspace{4mm}
\noindent
with probability
1 (further w. p. 1) for all $t\in[0, \bar T],$ where
${\rm exp}(At),$ $t\ge 0$ is the semigroup generated by the operator $A$.

Consider eigenvalues $\lambda_i$ and 
eigenfunctions $e_i(x)$ of the covariance operator $Q,$ where
$i=(i_1,\ldots,i_d)$ $\in$ $J,$
$x=(x_1,\ldots,x_d),$
and 
$J=\{i:\ i\in \mathbb{N}^d,\ \hbox{and}\ \lambda_i>0\}$.

The series representation of 
the $Q$-Wiener process has the following form \cite{7}

\vspace{1mm}
$$
{\bf W}_t=\sum\limits_{i\in J} e_i\sqrt{\lambda_i}
{\bf w}_t^{(i)}\ \ \ \hbox{or}\ \ \ {\bf W}_t=\sum\limits_{i\in J_M} 
e_i\ \langle e_i,{\bf W}_t\rangle_U,
$$

\vspace{3mm}
\noindent
where $t\in[0, \bar{T}],$\ 
${\bf w}_t^{(i)}$ $(i\in J)$ are independent standard
Wiener processes and $\langle \cdot,\cdot\rangle_U$ is a 
scalar product in $U.$

Note that eigenfunctions $e_i,$ $i\in J$ form
an orthonormal basis of $U$ \cite{7}.
Consider the finite-di\-men\-si\-onal approximation of ${\bf W}_t$ \cite{7}

\begin{equation}
\label{yy.1}
{\bf W}^M_t=\sum\limits_{i\in J_M} e_i\sqrt{\lambda_i}
{\bf w}_t^{(i)},\ \ \ t\in[0, \bar{T}],
\end{equation}

\vspace{2mm}
\noindent
where 
\begin{equation}
\label{ggg}
J_M=\{i:\ 1\le i_1,\ldots,i_d\le M,\ \hbox{and}\ 
\lambda_i>0\}.
\end{equation} 

\vspace{4mm}

{\bf Remark 1.}\ {\it Obviously, without the loss of generality
we can write $J_M=\{1, 2,\ldots, M\}.$}

Let $\Delta>0,$ $\tau_p=p\Delta$ $(p=0, 1,\ldots, N),$ 
and $N\Delta=\bar T.$ Consider the following 
exponential Milstein and Wagner--Platen numerical schemes 
for the SPDE (\ref{xx1}) \cite{2}, \cite{3} 

$$
Y_{p+1}={\rm exp}\left(A\Delta\right)\left(
Y_p+\Delta F(Y_p)+\int\limits_{\tau_p}^{\tau_{p+1}}B(Y_p)d{\bf W}_s
+\right.
$$

\begin{equation}
\label{fff100}
\left.+\int\limits_{\tau_p}^{\tau_{p+1}}B'(Y_p)
\left(\int\limits_{\tau_p}^{s}B(Y_p)d{\bf W}_{\tau}\right)
d{\bf W}_s\right),
\end{equation}

\vspace{5mm}
\centerline{(an exponential Milstein scheme)}

\vspace{9mm}

$$
Y_{p+1}={\rm exp}\left(\frac{A\Delta}{2}\right)\times
$$

\vspace{1mm}
$$
\times \left(
{\rm exp}\left(\frac{A\Delta}{2}\right)Y_p+\Delta F(Y_p)
+
\frac{\Delta^2}{2}
F'(Y_p)\biggl(AY_p+F(Y_p)\biggr)+
\int\limits_{\tau_p}^{\tau_{p+1}}B(Y_p)d{\bf W}_s
+\right.
$$

\vspace{1mm}
$$
+
\frac{\Delta^2}{4}
\sum\limits_{i\in J}\lambda_i F''(Y_p)\biggl(B(Y_p)e_i,B(Y_p)e_i\biggr)+
\int\limits_{\tau_p}^{\tau_{p+1}}B'(Y_p)
\left(\int\limits_{\tau_p}^{s}B(Y_p)d{\bf W}_{\tau}\right)
d{\bf W}_s+
$$

\vspace{1mm}
$$
+
A\left(\int\limits_{\tau_p}^{\tau_{p+1}}
\int\limits_{\tau_p}^{s}
B(Y_p)d{\bf W}_{\tau}ds-\frac{\Delta}{2}
\int\limits_{\tau_p}^{\tau_{p+1}}B(Y_p)d{\bf W}_s\right)+
$$

\vspace{3mm}
$$
+\Delta
\int\limits_{\tau_p}^{\tau_{p+1}}B'(Y_p)
\biggl(AY_p+F(Y_p)\biggr)d{\bf W}_s-
\int\limits_{\tau_p}^{\tau_{p+1}}\int\limits_{\tau_p}^{s}B'(Y_p)
\biggl(AY_p+F(Y_p)\biggr)d{\bf W}_{\tau}ds+
$$

\vspace{1mm}
$$
+
\frac{1}{2}\int\limits_{\tau_p}^{\tau_{p+1}}B''(Y_p)
\left(\int\limits_{\tau_p}^{s}B(Y_p)d{\bf W}_{\tau},
\int\limits_{\tau_p}^{s}B(Y_p)d{\bf W}_{\tau}
\right)
d{\bf W}_s+
$$

\vspace{1mm}
$$
+\int\limits_{\tau_p}^{\tau_{p+1}}
F'(Y_p)\left(\int\limits_{\tau_p}^s
B(Y_p)d{\bf W}_{\tau}\right)ds+
$$

\vspace{1mm}
\begin{equation}
\label{fff}
\left.+
\int\limits_{\tau_p}^{\tau_{p+1}}B'(Y_p)
\left(\int\limits_{\tau_p}^{s}B'(Y_p)
\left(\int\limits_{\tau_p}^{\tau}B(Y_p)
d{\bf W}_{\theta}\right)
d{\bf W}_{\tau}
\right)
d{\bf W}_s\right),
\end{equation}

\vspace{5mm}
\centerline{(an exponential Wagner--Platen scheme)}

\vspace{9mm}
\noindent 
where $Y_p$ is an approximation of $X_{\tau_p}$ (mild solution (\ref{mild})
at the time moment $\tau_p$), $p=0, 1,\ldots,N,$ and
$B'$, $B'',$ $F'$, $F''$ are
Fr\^{e}chet derivatives. 
In addition to the temporal discretization, 
the implementation of numerical schemes 
(\ref{fff100}) and (\ref{fff}) also requires a 
finite-dimensional approximation of the spaces $H$, $U$. 
Further, we will consider this approximation 
only for the space $U$.

Let us consider the following iterated It\^{o} stochastic integrals

$$
I_{(1)T,t}^{(r_1)}=\int\limits_t^Td{\bf w}_{t_1}^{(r_1)},
$$

$$
I_{(10)T,t}^{(r_1 0)}=
\int\limits_t^T\int\limits_t^{t_2} d{\bf w}_{t_1}^{(r_1)}
dt_2,\ \ \
I_{(01)T,t}^{(0 r_2)}=\int\limits_t^T\int\limits_t^{t_2} dt_1
d{\bf w}_{t_2}^{(r_2)},
$$

$$
I_{(11)T,t}^{(r_1 r_2)}=
\int\limits_t^T\int\limits_t^{t_2} d{\bf w}_{t_1}^{(r_1)}
d{\bf w}_{t_2}^{(r_2)},
$$

\vspace{2mm}

$$
I_{(111)T,t}^{(r_1 r_2 r_3)}=
\int\limits_t^T\int\limits_t^{t_3}\int\limits_t^{t_2} 
d{\bf w}_{t_1}^{(r_1)}
d{\bf w}_{t_2}^{(r_2)}
d{\bf w}_{t_3}^{(r_3)},
$$

\vspace{5mm}
\noindent
where $r_1, r_2, r_3\in J_M,$\ \ $0\le t<T\le \bar T,$
and $J_M$ is defined by (\ref{ggg}).

Let us replace the infinite-dimensional
$Q$-Wiener process in the iterated stochastic 
integrals from (\ref{fff100}), 
(\ref{fff}) by its finite-dimensional approximation (\ref{yy.1}).
Then we have w. p. 1

\begin{equation}
\label{ggg0}
\int\limits_{\tau_p}^{\tau_{p+1}}B(Y_p)d{\bf W}_s^M=
\sum\limits_{r_1\in J_M}
B(Y_p)e_{r_1}\sqrt{\lambda_{r_1}}
I_{(1)\tau_{p+1},\tau_p}^{(r_1)},
\end{equation}

\vspace{5mm}

$$
A\left(\int\limits_{\tau_p}^{\tau_{p+1}}
\int\limits_{\tau_p}^{s}
B(Y_p)d{\bf W}_{\tau}^M ds-\frac{\Delta}{2}
\int\limits_{\tau_p}^{\tau_{p+1}}B(Y_p)d{\bf W}_s^M\right)=
A\int\limits_{\tau_p}^{\tau_{p+1}}
B(Y_p)\left(\frac{\tau_{p+1}}{2}-s+\frac{\tau_p}{2}\right)
d{\bf W}_{s}^M=
$$

\vspace{1mm}
\begin{equation}
\label{ggg1}
=\sum\limits_{r_1\in J_M}
AB(Y_p)e_{r_1}\sqrt{\lambda_{r_1}}
\left(\frac{\Delta}{2}I_{(1)\tau_{p+1},\tau_p}^{(r_1)}-
I_{(01)\tau_{p+1},\tau_p}^{(0 r_1)}\right),
\end{equation}

\vspace{7mm}

$$
\Delta
\int\limits_{\tau_p}^{\tau_{p+1}}B'(Y_p)
\biggl(AY_p+F(Y_p)\biggr)d{\bf W}_s^M-
\int\limits_{\tau_p}^{\tau_{p+1}}\int\limits_{\tau_p}^{s}B'(Y_p)
\biggl(AY_p+F(Y_p)\biggr)d{\bf W}_{\tau}^Mds=
$$

\vspace{1mm}
$$
=\int\limits_{\tau_p}^{\tau_{p+1}}B'(Y_p)
\int\limits_{\tau_p}^{s}
\biggl(AY_p+F(Y_p)\biggr)d\tau d{\bf W}_{s}^M=
$$

\vspace{1mm}
\begin{equation}
\label{ggg2}
=
\sum\limits_{r_1\in J_M}
B'(Y_p)\biggl(AY_p+F(Y_p)\biggr)
e_{r_1}\sqrt{\lambda_{r_1}}I_{(01)\tau_{p+1},\tau_p}^{(0 r_1)},
\end{equation}

\vspace{7mm}

$$
\int\limits_{\tau_p}^{\tau_{p+1}}
F'(Y_p)\left(\int\limits_{\tau_p}^s
B(Y_p)d{\bf W}_{\tau}^M\right)ds=
$$

\vspace{1mm}
\begin{equation}
\label{ggg3}
=\sum\limits_{r_1\in J_M}
F'(Y_p)B(Y_p)e_{r_1}\sqrt{\lambda_{r_1}}
\left(\Delta I_{(1)\tau_{p+1},\tau_p}^{(r_1)}-
I_{(01)\tau_{p+1},\tau_p}^{(0 r_1)}\right),
\end{equation}

\vspace{7mm}

$$
\int\limits_{\tau_p}^{\tau_{p+1}}
B'(Y_p)\left(\int\limits_{\tau_p}^s
B(Y_p)d{\bf W}_{\tau}^M\right)d{\bf W}_s^M=
$$

\vspace{1mm}
\begin{equation}
\label{ch100}
=\sum\limits_{r_1,r_2\in J_M}
B'(Y_p)\left(B(Y_p)e_{r_1}\right)
e_{r_2}\sqrt{\lambda_{r_1}\lambda_{r_2}}I_{(11)\tau_{p+1},\tau_p}^{(r_1 r_2)},
\end{equation}

\vspace{9mm}

$$
\int\limits_{\tau_p}^{\tau_{p+1}}
B'(Y_p)\left(\int\limits_{\tau_p}^s
B'(Y_p)
\left(\int\limits_{\tau_p}^{\tau}
B(Y_p)
d{\bf W}_{\theta}^M\right)d{\bf W}_{\tau}^M\right)d{\bf W}_{s}^M=
$$

\vspace{1mm}
\begin{equation}
\label{ch101}
=\sum\limits_{r_1,r_2,r_3\in J_M}
B'(Y_p)\left(B'(Y_p)\left(B(Y_p)e_{r_1}\right)
e_{r_2}\right)e_{r_3}
\sqrt{\lambda_{r_1}\lambda_{r_2}\lambda_{r_3}}
I_{(111)\tau_{p+1},\tau_p}^{(r_1 r_2 r_3)},
\end{equation}

\vspace{7mm}

$$
\int\limits_{\tau_p}^{\tau_{p+1}}B''(Y_p)
\left(\int\limits_{\tau_p}^{s}B(Y_p)d{\bf W}_{\tau}^M,
\int\limits_{\tau_p}^{s}B(Y_p)d{\bf W}_{\tau}^M
\right)
d{\bf W}_s^M=
$$

\vspace{3mm}
$$
=
\sum_{r_1,r_2,r_3\in J_M}B''(Y_p)
\left(B(Y_p)e_{r_1}, B(Y_p)e_{r_2}\right)e_{r_3}
\sqrt{\lambda_{r_1}\lambda_{r_2}\lambda_{r_3}}\times
$$

\vspace{1mm}
\begin{equation}
\label{ggg4}
\times
\int\limits_{\tau_p}^{\tau_{p+1}}
\left(\int\limits_{\tau_p}^s d{\bf w}_{\tau}^{(r_1)}   
\int\limits_{\tau_p}^s d{\bf w}_{\tau}^{(r_2)}\right)
d{\bf w}_s^{(r_3)}.
\end{equation}

\vspace{7mm}

Note that in (\ref{ggg1})--(\ref{ggg3}) we used the It\^{o} formula.
Moreover, using the It\^{o} formula we obtain

\vspace{1mm}
\begin{equation}
\label{da12}
\int\limits_{\tau_p}^s d{\bf w}_{\tau}^{(r_1)}   
\int\limits_{\tau_p}^s d{\bf w}_{\tau}^{(r_2)}
=I_{(11)s,\tau_p}^{(r_1 r_2)}+I_{(11)s,\tau_p}^{(r_2 r_1)}
+
{\bf 1}_{\{r_1=r_2\}}(s-\tau_p)\ \ \ \hbox{w.\ p.\ 1.}
\end{equation}

\vspace{3mm}

From (\ref{da12}) we have

\vspace{1mm}

\begin{equation}
\label{jjj2}
\int\limits_{\tau_p}^{\tau_{p+1}}
\left(\int\limits_{\tau_p}^s d{\bf w}_{\tau}^{(r_1)}   
\int\limits_{\tau_p}^s d{\bf w}_{\tau}^{(r_2)}\right)
d{\bf w}_s^{(r_3)}=I_{(111)\tau_{p+1},\tau_p}^{(r_1 r_2 r_3)}
+I_{(111)\tau_{p+1},\tau_p}^{(r_2 r_1 r_3)}+
{\bf 1}_{\{r_1=r_2\}}I_{(01)\tau_{p+1},\tau_p}^{(0 r_3)}
\ \ \ \hbox{w.\ p.\ 1.}
\end{equation}

\vspace{5mm}

After substituting (\ref{jjj2}) into (\ref{ggg4}), we obtain

\vspace{2mm}
$$
\int\limits_{\tau_p}^{\tau_{p+1}}B''(Y_p)
\left(\int\limits_{\tau_p}^{s}B(Y_p)d{\bf W}_{\tau}^M,
\int\limits_{\tau_p}^{s}B(Y_p)d{\bf W}_{\tau}^M
\right)
d{\bf W}_s^M=
$$

\vspace{1mm}
$$
=
\sum_{r_1,r_2,r_3\in J_M}B''(Z)\left(B(Z)e_{r_1}, B(Z)e_{r_2}\right)e_{r_3}
\sqrt{\lambda_{r_1}\lambda_{r_2}\lambda_{r_3}}\times
$$

\vspace{1mm}
\begin{equation}
\label{da1}
\times \left(I_{(111)\tau_{p+1},\tau_p}^{(r_1 r_2 r_3)}
+I_{(111)\tau_{p+1},\tau_p}^{(r_2 r_1 r_3)}+
{\bf 1}_{\{r_1=r_2\}}
I_{(01)\tau_{p+1},\tau_p}^{(0 r_3)}\right)\ \ \ \hbox{w.\ p.\ 1.}
\end{equation}

\vspace{7mm}

Thus, for the implementation of numerical schemes 
(\ref{fff100}) and (\ref{fff})
we need to approximate the following collection of iterated 
It\^{o} stochastic integrals

$$
I_{(1)T,t}^{(r_1)},\ \ \
I_{(01)T,t}^{(0 r_1)},\ \ \ 
I_{(11)T,t}^{(r_1 r_2)},\ \ \
I_{(111)T,t}^{(r_1 r_2 r_3)},
$$

\vspace{3mm}
\noindent
where $r_1, r_2, r_3\in J_M,$\ \  $0\le t<T\le \bar T.$

The monographs
\cite{24} (Chapters 5, 6) and \cite{24a} or \cite{24aa} (Chapters 1, 2, and 5)
(also see \cite{11}-\cite{23}, \cite{21}-\cite{arxiv-25})
are devoted to the constructing of efficient methods (based on generalized multiple Fourier series)
of the mean-square approximation of iterated It\^{o} stochastic 
integrals with respect to components of the finite-dimensional
Wiener process. 
These results are also adapted for iterated Stratonovich stochastic integrals
\cite{16}-\cite{24aa}, \cite{30}, \cite{arxiv-2}, \cite{arxiv-4}-\cite{arxiv-11},
\cite{arxiv-15}, \cite{arxiv-18}, \cite{arxiv-19}, \cite{arxiv-23},
\cite{arxiv-25}.
In Sect.~3, we consider a very short review of the results from 
monographs \cite{24} (Chapters 5, 6) and \cite{24a} or \cite{24aa} (Chapters 1, 2, and 5).

\vspace{5mm}

\section{Method of Approximation of Iterated 
It\^{o} Stohastic integrals Based on Generalized Multiple Fourier Series. The 
Case of Multiple Fo\-u\-ri\-er--Le-\\ 
gendre Series}

\vspace{5mm}

Consider the following 
iterated It\^{o} stochastic integrals

\begin{equation}
\label{sodom20}
J[\psi^{(k)}]_{T,t}^{(i_1\ldots i_k)}
=\int\limits_t^T\psi_k(t_k) \ldots \int\limits_t^{t_{2}}
\psi_1(t_1) d{\bf w}_{t_1}^{(i_1)}\ldots
d{\bf w}_{t_k}^{(i_k)},
\end{equation}

\vspace{2mm}
\noindent
where $0\le t<T \le \bar{T},$ every $\psi_l(\tau)\ (l=1,\ldots,k)$ is 
a continuous non-random function on $[t, T]$,
${\bf w}_{\tau}^{(i)}$ ($i=1,\ldots,m)$ are independent
standard Wiener processes, 
${\bf w}_{\tau}^{(0)}=\tau,$\
$i_1,\ldots,i_k=0, 1,\ldots,m.$ 
The case $\psi_1(\tau),\ldots,\psi_k(\tau)\in L_2([t, T])$
will be considered in Theorem~2 (see below).

Suppose that $\{\phi_j(x)\}_{j=0}^{\infty}$
is a complete orthonormal system of functions in the space
$L_2([t, T])$. 
Define the following function on the hypercube $[t, T]^k$

\vspace{-1mm}
\begin{equation}
\label{ppp}
K(t_1,\ldots,t_k)=
\begin{cases}
\psi_1(t_1)\ldots \psi_k(t_k),\ t_1<\ldots<t_k\\
~\\
~\\
0,\ \hbox{\rm otherwise}
\end{cases}
=\ \
\prod\limits_{l=1}^k
\psi_l(t_l)\ \prod\limits_{l=1}^{k-1}{\bf 1}_{\{t_l<t_{l+1}\}},\ 
\end{equation}

\vspace{4mm}
\noindent
where $t_1,\ldots,t_k\in [t, T]$ for $k\ge 2$ and 
$K(t_1)\equiv\psi_1(t_1)$ for $t_1\in[t, T].$ Here 
${\bf 1}_A$ is the indicator of the set $A$.

The function $K(t_1,\ldots,t_k)$ is piecewise continuous on the 
hypercube $[t, T]^k.$
At this situation it is well known that the generalized 
multiple Fourier series 
of $K(t_1,\ldots,t_k)\in L_2([t, T]^k)$ converges
to $K(t_1,\ldots,t_k)$ on the hypercube $[t, T]^k$ in 
the mean-square sense, i.e.

\begin{equation}
\label{sos1z}
\hbox{\vtop{\offinterlineskip\halign{
\hfil#\hfil\cr
{\rm lim}\cr
$\stackrel{}{{}_{p_1,\ldots,p_k\to \infty}}$\cr
}} }\Biggl\Vert
K(t_1,\ldots,t_k)-
\sum_{j_1=0}^{p_1}\ldots \sum_{j_k=0}^{p_k}
C_{j_k\ldots j_1}\prod_{l=1}^{k} \phi_{j_l}(t_l)\Biggr
\Vert_{L_2([t, T]^k)}=0,
\end{equation}

\vspace{2mm}
\noindent
where

\vspace{-3mm}
\begin{equation}
\label{ppppa}
C_{j_k\ldots j_1}=\int\limits_{[t,T]^k}
K(t_1,\ldots,t_k)\prod_{l=1}^{k}\phi_{j_l}(t_l)dt_1\ldots dt_k
\end{equation}

\vspace{5mm}
\noindent
is the Fourier coefficient and

$$
\left\Vert f\right\Vert_{L_2([t, T]^k)}=\left(\int\limits_{[t,T]^k}
f^2(t_1,\ldots,t_k)dt_1\ldots dt_k\right)^{1/2}.
$$

\vspace{4mm}

Consider the discretization $\{\tau_j\}_{j=0}^N$ of $[t,T]$ such that

\begin{equation}
\label{1111}
t=\tau_0<\ldots <\tau_N=T,\ \ \ \
\Delta_N=
\hbox{\vtop{\offinterlineskip\halign{
\hfil#\hfil\cr
{\rm max}\cr
$\stackrel{}{{}_{0\le j\le N-1}}$\cr
}} }\Delta\tau_j\to 0\ \ \hbox{if}\ \ N\to \infty,\ \ \ \
\Delta\tau_j=\tau_{j+1}-\tau_j.
\end{equation}

\vspace{5mm}

{\bf Theorem 1} \cite{11} (2006) (also see \cite{9}, \cite{arxiv-20},
\cite{12}-\cite{arxiv-25}). {\it 
Suppose that
every $\psi_l(\tau)$ $(l=1,\ldots, k)$ is a conti\-nu\-ous 
non-random function on $[t, T]$ and
$\{\phi_j(x)\}_{j=0}^{\infty}$ is a complete orthonormal system  
of continuous functions in $L_2([t,T]).$ Then

\vspace{2mm}
$$
J[\psi^{(k)}]_{T,t}^{(i_1\ldots i_k)}=
\hbox{\vtop{\offinterlineskip\halign{
\hfil#\hfil\cr
{\rm l.i.m.}\cr
$\stackrel{}{{}_{p_1,\ldots,p_k\to \infty}}$\cr
}} }\sum_{j_1=0}^{p_1}\ldots\sum_{j_k=0}^{p_k}
C_{j_k\ldots j_1}\Biggl(
\prod_{l=1}^k\zeta_{j_l}^{(i_l)}-
\Biggr.
$$

\vspace{3mm}
\begin{equation}
\label{tyyy}
-\ \Biggl.
\hbox{\vtop{\offinterlineskip\halign{
\hfil#\hfil\cr
{\rm l.i.m.}\cr
$\stackrel{}{{}_{N\to \infty}}$\cr
}} }\sum_{(l_1,\ldots,l_k)\in {\rm G}_k}
\phi_{j_{1}}(\tau_{l_1})
\Delta{\bf w}_{\tau_{l_1}}^{(i_1)}\ldots
\phi_{j_{k}}(\tau_{l_k})
\Delta{\bf w}_{\tau_{l_k}}^{(i_k)}\Biggr),
\end{equation}

\vspace{5mm}
\noindent
where

$$
{\rm G}_k={\rm H}_k\backslash{\rm L}_k,\ \ \
{\rm H}_k=\biggl\{(l_1,\ldots,l_k):\ l_1,\ldots,l_k=0,\ 1,\ldots,N-1\biggr\},
$$

$$
{\rm L}_k=\biggl\{(l_1,\ldots,l_k):\ l_1,\ldots,l_k=0,\ 1,\ldots,N-1;\
l_g\ne l_r\ (g\ne r);\ g, r=1,\ldots,k\biggr\},
$$

\vspace{5mm}
\noindent
${\rm l.i.m.}$ is a limit in the mean-square sense{\rm ,} 
$i_1,\ldots,i_k=0,1,\ldots,m,$

\begin{equation}
\label{rr23}
\zeta_{j}^{(i)}=
\int\limits_t^T \phi_{j}(s) d{\bf w}_s^{(i)}
\end{equation} 

\vspace{3mm}
\noindent
are independent standard Gaussian random variables
for various
$i$ or $j$ {\rm(}if $i\ne 0${\rm),}
$C_{j_k\ldots j_1}$ is the Fourier coefficient {\rm(\ref{ppppa}),}
$\Delta{\bf w}_{\tau_{j}}^{(i)}=
{\bf w}_{\tau_{j+1}}^{(i)}-{\bf w}_{\tau_{j}}^{(i)}$
$(i=0,\ 1,\ldots,m),$\
$\left\{\tau_{j}\right\}_{j=0}^{N}$ 
is the discretization of
$[t,T],$ which satisfies the condition {\rm (\ref{1111})}.
}

Note that in \cite{11}-\cite{24aa}, \cite{arxiv-1}
the version of Theorem 1 for
systems of Haar and Rademacher--Walsh 
functions has been considered. 
Another version of Theorem 1 related to the application of 
complete orthonormal with weight $r (t_1)\ldots r(t_k)\ge 0$ systems of 
functions in $L_2 ([t, T]^k)$ has been 
considered in \cite{24}-\cite{24aa}, \cite{arxiv-13}.
A generalization of Theorem~1 to the case
of an arbitrary complete orthonormal system  
of functions $\{\phi_j(x)\}_{j=0}^{\infty}$ in the space $L_2([t,T])$ as well as 
$\psi_1(\tau),\ldots,\psi_k(\tau)\in L_2([t, T])$ will be considered below (see Theorem~2).

In order to evaluate the significance of Theorem 1 for practice we will
demonstrate its transformed particular cases for 
$k=1,\ldots,5$ \cite{11}-\cite{arxiv-25}
(cases $k=6,$ $7$ and $k>7$
can be found in \cite{12}-\cite{arxiv-25})

\begin{equation}
\label{a1}
J[\psi^{(1)}]_{T,t}^{(i_1)}
=\hbox{\vtop{\offinterlineskip\halign{
\hfil#\hfil\cr
{\rm l.i.m.}\cr
$\stackrel{}{{}_{p_1\to \infty}}$\cr
}} }\sum_{j_1=0}^{p_1}
C_{j_1}\zeta_{j_1}^{(i_1)},
\end{equation}

\vspace{2mm}

\begin{equation}
\label{a2}
J[\psi^{(2)}]_{T,t}^{(i_1 i_2)}
=\hbox{\vtop{\offinterlineskip\halign{
\hfil#\hfil\cr
{\rm l.i.m.}\cr
$\stackrel{}{{}_{p_1,p_2\to \infty}}$\cr
}} }\sum_{j_1=0}^{p_1}\sum_{j_2=0}^{p_2}
C_{j_2j_1}\Biggl(\zeta_{j_1}^{(i_1)}\zeta_{j_2}^{(i_2)}
-{\bf 1}_{\{i_1=i_2\ne 0\}}
{\bf 1}_{\{j_1=j_2\}}\Biggr),
\end{equation}

\vspace{5mm}
$$
J[\psi^{(3)}]_{T,t}^{(i_1 i_2 i_3)}=
\hbox{\vtop{\offinterlineskip\halign{
\hfil#\hfil\cr
{\rm l.i.m.}\cr
$\stackrel{}{{}_{p_1,p_2,p_3\to \infty}}$\cr
}} }\sum_{j_1=0}^{p_1}\sum_{j_2=0}^{p_2}\sum_{j_3=0}^{p_3}
C_{j_3j_2j_1}\Biggl(
\zeta_{j_1}^{(i_1)}\zeta_{j_2}^{(i_2)}\zeta_{j_3}^{(i_3)}
-\Biggr.
$$
\begin{equation}
\label{a3}
-{\bf 1}_{\{i_1=i_2\ne 0\}}
{\bf 1}_{\{j_1=j_2\}}
\zeta_{j_3}^{(i_3)}
-{\bf 1}_{\{i_2=i_3\ne 0\}}
{\bf 1}_{\{j_2=j_3\}}
\zeta_{j_1}^{(i_1)}
\Biggl.-{\bf 1}_{\{i_1=i_3\ne 0\}}
{\bf 1}_{\{j_1=j_3\}}
\zeta_{j_2}^{(i_2)}\Biggr),
\end{equation}

\vspace{7mm}

$$
J[\psi^{(4)}]_{T,t}^{(i_1\ldots i_4)}
=
\hbox{\vtop{\offinterlineskip\halign{
\hfil#\hfil\cr
{\rm l.i.m.}\cr
$\stackrel{}{{}_{p_1,\ldots,p_4\to \infty}}$\cr
}} }\sum_{j_1=0}^{p_1}\ldots\sum_{j_4=0}^{p_4}
C_{j_4\ldots j_1}\Biggl(
\prod_{l=1}^4\zeta_{j_l}^{(i_l)}
\Biggr.
-
$$
$$
-
{\bf 1}_{\{i_1=i_2\ne 0\}}
{\bf 1}_{\{j_1=j_2\}}
\zeta_{j_3}^{(i_3)}
\zeta_{j_4}^{(i_4)}
-
{\bf 1}_{\{i_1=i_3\ne 0\}}
{\bf 1}_{\{j_1=j_3\}}
\zeta_{j_2}^{(i_2)}
\zeta_{j_4}^{(i_4)}-
$$
$$
-
{\bf 1}_{\{i_1=i_4\ne 0\}}
{\bf 1}_{\{j_1=j_4\}}
\zeta_{j_2}^{(i_2)}
\zeta_{j_3}^{(i_3)}
-
{\bf 1}_{\{i_2=i_3\ne 0\}}
{\bf 1}_{\{j_2=j_3\}}
\zeta_{j_1}^{(i_1)}
\zeta_{j_4}^{(i_4)}-
$$
$$
-
{\bf 1}_{\{i_2=i_4\ne 0\}}
{\bf 1}_{\{j_2=j_4\}}
\zeta_{j_1}^{(i_1)}
\zeta_{j_3}^{(i_3)}
-
{\bf 1}_{\{i_3=i_4\ne 0\}}
{\bf 1}_{\{j_3=j_4\}}
\zeta_{j_1}^{(i_1)}
\zeta_{j_2}^{(i_2)}+
$$

\vspace{-3mm}
$$
+
{\bf 1}_{\{i_1=i_2\ne 0\}}
{\bf 1}_{\{j_1=j_2\}}
{\bf 1}_{\{i_3=i_4\ne 0\}}
{\bf 1}_{\{j_3=j_4\}}
+
$$

\vspace{-3mm}
$$
+
{\bf 1}_{\{i_1=i_3\ne 0\}}
{\bf 1}_{\{j_1=j_3\}}
{\bf 1}_{\{i_2=i_4\ne 0\}}
{\bf 1}_{\{j_2=j_4\}}+
$$
\begin{equation}
\label{a4}
+\Biggl.
{\bf 1}_{\{i_1=i_4\ne 0\}}
{\bf 1}_{\{j_1=j_4\}}
{\bf 1}_{\{i_2=i_3\ne 0\}}
{\bf 1}_{\{j_2=j_3\}}\Biggr),
\end{equation}

\vspace{7mm}

$$
J[\psi^{(5)}]_{T,t}^{(i_1\ldots i_5)}
=\hbox{\vtop{\offinterlineskip\halign{
\hfil#\hfil\cr
{\rm l.i.m.}\cr
$\stackrel{}{{}_{p_1,\ldots,p_5\to \infty}}$\cr
}} }\sum_{j_1=0}^{p_1}\ldots\sum_{j_5=0}^{p_5}
C_{j_5\ldots j_1}\Biggl(
\prod_{l=1}^5\zeta_{j_l}^{(i_l)}
-\Biggr.
$$
$$
-
{\bf 1}_{\{i_1=i_2\ne 0\}}
{\bf 1}_{\{j_1=j_2\}}
\zeta_{j_3}^{(i_3)}
\zeta_{j_4}^{(i_4)}
\zeta_{j_5}^{(i_5)}-
{\bf 1}_{\{i_1=i_3\ne 0\}}
{\bf 1}_{\{j_1=j_3\}}
\zeta_{j_2}^{(i_2)}
\zeta_{j_4}^{(i_4)}
\zeta_{j_5}^{(i_5)}-
$$
$$
-
{\bf 1}_{\{i_1=i_4\ne 0\}}
{\bf 1}_{\{j_1=j_4\}}
\zeta_{j_2}^{(i_2)}
\zeta_{j_3}^{(i_3)}
\zeta_{j_5}^{(i_5)}-
{\bf 1}_{\{i_1=i_5\ne 0\}}
{\bf 1}_{\{j_1=j_5\}}
\zeta_{j_2}^{(i_2)}
\zeta_{j_3}^{(i_3)}
\zeta_{j_4}^{(i_4)}-
$$
$$
-
{\bf 1}_{\{i_2=i_3\ne 0\}}
{\bf 1}_{\{j_2=j_3\}}
\zeta_{j_1}^{(i_1)}
\zeta_{j_4}^{(i_4)}
\zeta_{j_5}^{(i_5)}-
{\bf 1}_{\{i_2=i_4\ne 0\}}
{\bf 1}_{\{j_2=j_4\}}
\zeta_{j_1}^{(i_1)}
\zeta_{j_3}^{(i_3)}
\zeta_{j_5}^{(i_5)}-
$$
$$
-
{\bf 1}_{\{i_2=i_5\ne 0\}}
{\bf 1}_{\{j_2=j_5\}}
\zeta_{j_1}^{(i_1)}
\zeta_{j_3}^{(i_3)}
\zeta_{j_4}^{(i_4)}
-{\bf 1}_{\{i_3=i_4\ne 0\}}
{\bf 1}_{\{j_3=j_4\}}
\zeta_{j_1}^{(i_1)}
\zeta_{j_2}^{(i_2)}
\zeta_{j_5}^{(i_5)}-
$$
$$
-
{\bf 1}_{\{i_3=i_5\ne 0\}}
{\bf 1}_{\{j_3=j_5\}}
\zeta_{j_1}^{(i_1)}
\zeta_{j_2}^{(i_2)}
\zeta_{j_4}^{(i_4)}
-{\bf 1}_{\{i_4=i_5\ne 0\}}
{\bf 1}_{\{j_4=j_5\}}
\zeta_{j_1}^{(i_1)}
\zeta_{j_2}^{(i_2)}
\zeta_{j_3}^{(i_3)}+
$$
$$
+
{\bf 1}_{\{i_1=i_2\ne 0\}}
{\bf 1}_{\{j_1=j_2\}}
{\bf 1}_{\{i_3=i_4\ne 0\}}
{\bf 1}_{\{j_3=j_4\}}\zeta_{j_5}^{(i_5)}+
{\bf 1}_{\{i_1=i_2\ne 0\}}
{\bf 1}_{\{j_1=j_2\}}
{\bf 1}_{\{i_3=i_5\ne 0\}}
{\bf 1}_{\{j_3=j_5\}}\zeta_{j_4}^{(i_4)}+
$$
$$
+
{\bf 1}_{\{i_1=i_2\ne 0\}}
{\bf 1}_{\{j_1=j_2\}}
{\bf 1}_{\{i_4=i_5\ne 0\}}
{\bf 1}_{\{j_4=j_5\}}\zeta_{j_3}^{(i_3)}+
{\bf 1}_{\{i_1=i_3\ne 0\}}
{\bf 1}_{\{j_1=j_3\}}
{\bf 1}_{\{i_2=i_4\ne 0\}}
{\bf 1}_{\{j_2=j_4\}}\zeta_{j_5}^{(i_5)}+
$$
$$
+
{\bf 1}_{\{i_1=i_3\ne 0\}}
{\bf 1}_{\{j_1=j_3\}}
{\bf 1}_{\{i_2=i_5\ne 0\}}
{\bf 1}_{\{j_2=j_5\}}\zeta_{j_4}^{(i_4)}+
{\bf 1}_{\{i_1=i_3\ne 0\}}
{\bf 1}_{\{j_1=j_3\}}
{\bf 1}_{\{i_4=i_5\ne 0\}}
{\bf 1}_{\{j_4=j_5\}}\zeta_{j_2}^{(i_2)}+
$$
$$
+
{\bf 1}_{\{i_1=i_4\ne 0\}}
{\bf 1}_{\{j_1=j_4\}}
{\bf 1}_{\{i_2=i_3\ne 0\}}
{\bf 1}_{\{j_2=j_3\}}\zeta_{j_5}^{(i_5)}+
{\bf 1}_{\{i_1=i_4\ne 0\}}
{\bf 1}_{\{j_1=j_4\}}
{\bf 1}_{\{i_2=i_5\ne 0\}}
{\bf 1}_{\{j_2=j_5\}}\zeta_{j_3}^{(i_3)}+
$$
$$
+
{\bf 1}_{\{i_1=i_4\ne 0\}}
{\bf 1}_{\{j_1=j_4\}}
{\bf 1}_{\{i_3=i_5\ne 0\}}
{\bf 1}_{\{j_3=j_5\}}\zeta_{j_2}^{(i_2)}+
{\bf 1}_{\{i_1=i_5\ne 0\}}
{\bf 1}_{\{j_1=j_5\}}
{\bf 1}_{\{i_2=i_3\ne 0\}}
{\bf 1}_{\{j_2=j_3\}}\zeta_{j_4}^{(i_4)}+
$$
$$
+
{\bf 1}_{\{i_1=i_5\ne 0\}}
{\bf 1}_{\{j_1=j_5\}}
{\bf 1}_{\{i_2=i_4\ne 0\}}
{\bf 1}_{\{j_2=j_4\}}\zeta_{j_3}^{(i_3)}+
{\bf 1}_{\{i_1=i_5\ne 0\}}
{\bf 1}_{\{j_1=j_5\}}
{\bf 1}_{\{i_3=i_4\ne 0\}}
{\bf 1}_{\{j_3=j_4\}}\zeta_{j_2}^{(i_2)}+
$$
$$
+
{\bf 1}_{\{i_2=i_3\ne 0\}}
{\bf 1}_{\{j_2=j_3\}}
{\bf 1}_{\{i_4=i_5\ne 0\}}
{\bf 1}_{\{j_4=j_5\}}\zeta_{j_1}^{(i_1)}+
{\bf 1}_{\{i_2=i_4\ne 0\}}
{\bf 1}_{\{j_2=j_4\}}
{\bf 1}_{\{i_3=i_5\ne 0\}}
{\bf 1}_{\{j_3=j_5\}}\zeta_{j_1}^{(i_1)}+
$$
\begin{equation}
\label{a5}
+\Biggl.
{\bf 1}_{\{i_2=i_5\ne 0\}}
{\bf 1}_{\{j_2=j_5\}}
{\bf 1}_{\{i_3=i_4\ne 0\}}
{\bf 1}_{\{j_3=j_4\}}\zeta_{j_1}^{(i_1)}\Biggr),
\end{equation}

\vspace{7mm}
\noindent
where ${\bf 1}_A$ is the indicator of the set $A$.

Consider the generalization of (\ref{a1})--(\ref{a5}) 
for the case of an arbitrary $k$ $(k\in \mathbb{N})$ 
as well as for the case of an arbitrary complete orthonormal system 
of functions $\{\phi_j(x)\}_{j=0}^{\infty}$ in the space $L_2([t,T])$
and $\psi_1(\tau),\ldots,\psi_k(\tau) \in L_2([t, T])$.

In order to do this, let us
consider the unordered
set $\{1, 2, \ldots, k\}$ 
and separate it into two parts:
the first part consists of $r$ unordered 
pairs (sequence order of these pairs is also unimportant) and the 
second one consists of the 
remaining $k-2r$ numbers.
So, we have

\begin{equation}
\label{leto5007}
(\{
\underbrace{\{g_1, g_2\}, \ldots, 
\{g_{2r-1}, g_{2r}\}}_{\small{\hbox{part 1}}}
\},
\{\underbrace{q_1, \ldots, q_{k-2r}}_{\small{\hbox{part 2}}}
\}),
\end{equation}

\vspace{3mm}
\noindent
where 
$\{g_1, g_2, \ldots, 
g_{2r-1}, g_{2r}, q_1, \ldots, q_{k-2r}\}=\{1, 2, \ldots, k\},$
braces   
mean an unordered 
set and parentheses mean an ordered set.

We will say that (\ref{leto5007}) is a partition 
and consider the sum with respect to all possible
partitions

\vspace{1mm}
\begin{equation}
\label{leto5008}
\sum_{\stackrel{(\{\{g_1, g_2\}, \ldots, 
\{g_{2r-1}, g_{2r}\}\}, \{q_1, \ldots, q_{k-2r}\})}
{{}_{\{g_1, g_2, \ldots, 
g_{2r-1}, g_{2r}, q_1, \ldots, q_{k-2r}\}=\{1, 2, \ldots, k\}}}}
a_{g_1 g_2, \ldots, 
g_{2r-1} g_{2r}, q_1 \ldots q_{k-2r}}.
\end{equation}

\vspace{6mm}

{\bf Theorem 2}\ \cite{24a} (Sect.~1.11), \cite{arxiv-1} (Sect.~15). {\it Suppose that
$\{\phi_j(x)\}_{j=0}^{\infty}$ is an arbitrary complete orthonormal system  
of functions in the space $L_2([t,T])$ and 
$\psi_1(\tau),\ldots,\psi_k(\tau)\in L_2([t, T]).$ 
Then the following expansion

\vspace{2mm}

$$
J[\psi^{(k)}]_{T,t}^{(i_1\ldots i_k)}=
\hbox{\vtop{\offinterlineskip\halign{
\hfil#\hfil\cr
{\rm l.i.m.}\cr
$\stackrel{}{{}_{p_1,\ldots,p_k\to \infty}}$\cr
}} }
\sum\limits_{j_1=0}^{p_1}\ldots
\sum\limits_{j_k=0}^{p_k}
C_{j_k\ldots j_1}\Biggl(
\prod_{l=1}^k\zeta_{j_l}^{(i_l)}+\sum\limits_{r=1}^{[k/2]}
(-1)^r \times
\Biggr.
$$

\vspace{2mm}
\begin{equation}
\label{leto6000}
\times
\sum_{\stackrel{(\{\{g_1, g_2\}, \ldots, 
\{g_{2r-1}, g_{2r}\}\}, \{q_1, \ldots, q_{k-2r}\})}
{{}_{\{g_1, g_2, \ldots, 
g_{2r-1}, g_{2r}, q_1, \ldots, q_{k-2r}\}=\{1, 2, \ldots, k\}}}}
\prod\limits_{s=1}^r
{\bf 1}_{\{i_{g_{{}_{2s-1}}}=~i_{g_{{}_{2s}}}\ne 0\}}
\Biggl.{\bf 1}_{\{j_{g_{{}_{2s-1}}}=~j_{g_{{}_{2s}}}\}}
\prod_{l=1}^{k-2r}\zeta_{j_{q_l}}^{(i_{q_l})}\Biggr)
\end{equation}

\vspace{6mm}
\noindent
con\-verg\-ing in the mean-square sense is valid, where $[x]$ is an integer part of 
a real number $x;$ another notations are the same as in Theorem~{\rm 1}.}

\vspace{2mm}

In particular, from (\ref{leto6000}) for $k=5$ we obtain

\vspace{3mm}

$$
J[\psi^{(5)}]_{T,t}^{(i_1\ldots i_5)}=
\hbox{\vtop{\offinterlineskip\halign{
\hfil#\hfil\cr
{\rm l.i.m.}\cr
$\stackrel{}{{}_{p_1,\ldots,p_5\to \infty}}$\cr
}} }\sum_{j_1=0}^{p_1}\ldots\sum_{j_5=0}^{p_5}
C_{j_5\ldots j_1}\Biggl(
\prod_{l=1}^5\zeta_{j_l}^{(i_l)}-\Biggr.
$$

\vspace{3mm}

$$
-
\sum\limits_{\stackrel{(\{g_1, g_2\}, \{q_1, q_{2}, q_3\})}
{{}_{\{g_1, g_2, q_{1}, q_{2}, q_3\}=\{1, 2, 3, 4, 5\}}}}
{\bf 1}_{\{i_{g_{{}_{1}}}=~i_{g_{{}_{2}}}\ne 0\}}
{\bf 1}_{\{j_{g_{{}_{1}}}=~j_{g_{{}_{2}}}\}}
\prod_{l=1}^{3}\zeta_{j_{q_l}}^{(i_{q_l})}+
$$

\vspace{3mm}
$$
+
\sum_{\stackrel{(\{\{g_1, g_2\}, 
\{g_{3}, g_{4}\}\}, \{q_1\})}
{{}_{\{g_1, g_2, g_{3}, g_{4}, q_1\}=\{1, 2, 3, 4, 5\}}}}
{\bf 1}_{\{i_{g_{{}_{1}}}=~i_{g_{{}_{2}}}\ne 0\}}
{\bf 1}_{\{j_{g_{{}_{1}}}=~j_{g_{{}_{2}}}\}}
\Biggl.{\bf 1}_{\{i_{g_{{}_{3}}}=~i_{g_{{}_{4}}}\ne 0\}}
{\bf 1}_{\{j_{g_{{}_{3}}}=~j_{g_{{}_{4}}}\}}
\zeta_{j_{q_1}}^{(i_{q_1})}\Biggr).
$$

\vspace{8mm}

The last equality obviously agrees with
(\ref{a5}).
Note that the correctness of formulas (\ref{a1})--(\ref{a5}) 
can be 
verified 
by the fact that if 
$i_1=\ldots=i_5=i=1,\ldots,m$
and $\psi_1(s),\ldots,\psi_6(s)\equiv \psi(s)$,
then we can derive from (\ref{a1})--(\ref{a5}) the well known
equalities (the 
cases $k=2, 3$ were discussed in details in \cite{12}-\cite{24aa},
\cite{arxiv-1})

$$
J[\psi^{(1)}]_{T,t}^{(i)}
=\frac{1}{1!}\delta_{T,t}^{(i)},
$$

\vspace{3mm}
$$
J[\psi^{(2)}]_{T,t}^{(ii)}
=\frac{1}{2!}\left(\left(\delta^{(i)}_{T,t}\right)^2-\Delta_{T,t}\right),\
$$

\vspace{3mm}
$$
J[\psi^{(3)}]_{T,t}^{(iii)}
=\frac{1}{3!}\left(\left(\delta_{T,t}^{(i)}\right)^3-
3\delta_{T,t}^{(i)}\Delta_{T,t}\right),
$$

\vspace{3mm}
$$
J[\psi^{(4)}]_{T,t}^{(iiii)}
=\frac{1}{4!}\left(\left(\delta_{T,t}^{(i)}\right)^4-
6\left(\delta_{T,t}^{(i)}\right)^2\Delta_{T,t}
+3\Delta^2_{T,t}\right),\
$$

\vspace{4mm}
$$
J[\psi^{(5)}]_{T,t}^{(iiiii)}
=\frac{1}{5!}\left(\left(\delta_{T,t}^{(i)}\right)^5-
10\left(\delta_{T,t}^{(i)}\right)^3\Delta_{T,t}
+15\delta_{T,t}^{(i)}\Delta^2_{T,t}\right)
$$

\vspace{5mm}
\noindent
w. p. 1, where 

\vspace{-1mm}
$$
\delta_{T,t}^{(i)}=\int\limits_t^T\psi(s)d{\bf w}_s^{(i)},\ \ \
\Delta_{T,t}=\int\limits_t^T\psi^2(s)ds.
$$

\vspace{4mm}

The above equalities can be independently  
obtained using the It\^{o} formula and Hermite polynomials \cite{160}.

Assume that $J[\psi^{(k)}]_{T,t}^{(i_1\ldots i_k)p_1 \ldots p_k}$ 
is an approximation 
of the stochastic integral (\ref{sodom20}), 
which is the 
expression on the right-hand side of (\ref{leto6000}) before passing
to the limit $\hbox{\vtop{\offinterlineskip\halign{
\hfil#\hfil\cr
{\rm l.i.m.}\cr
$\stackrel{}{{}_{p_1,\ldots,p_k\to \infty}}$\cr
}} }.$\ \
Let us denote

\vspace{3mm}
$$
E^{(i_1\ldots i_k)p_1,\ldots,p_k}={\sf M}\left\{\biggl(
J[\psi^{(k)}]_{T,t}^{(i_1\ldots i_k)}-
J[\psi^{(k)}]_{T,t}^{(i_1\ldots i_k)p_1,\ldots,p_k}\biggr)^2\right\},
$$

\vspace{5mm}
\begin{equation}
\label{g123}
I_k=\Vert K\Vert^2_{L_2([t, T]^k)}=\int\limits_{[t,T]^k}
K^2(t_1,\ldots,t_k)dt_1\ldots dt_k,\ \ \ \
\left. E^{(i_1\ldots i_k)p_1,\ldots,p_k}\right|_{p_1
=\ldots=p_k=p}\stackrel{\sf def}{=}E^{(i_1\ldots i_k)p}.
\end{equation}  

\vspace{8mm}

In \cite{23}, \cite{24}-\cite{24aa}, 
\cite{arxiv-1}, \cite{arxiv-3}  it was shown that 

\vspace{2mm}
\begin{equation}
\label{star00011}
E_k^{(i_1\ldots i_k)p_1,\ldots,p_k}\le k!\left(I_k-\sum_{j_1=0}^{p_1}\ldots
\sum_{j_k=0}^{p_k}C^2_{j_k\ldots j_1}\right),
\end{equation}

\vspace{6mm}
\noindent
where $i_1,\ldots,i_k=1,\ldots,m$ for $0<T-t<\infty$ and
$i_1,\ldots,i_k=0, 1,\ldots,m$ for $0<T-t<1.$
Note that the
estimate (\ref{star00011}) is valid under the conditions of Theorem~2.

Let us consider some approximations
of iterated It\^{o} stochastic integrals
using Theorems 1, 2 and multiple Fourier--Legendre series.

The complete orthonormal system of Legendre polynomials in 
the space $L_2([t,T])$ looks as follows

\vspace{1mm}
\begin{equation}
\label{66}
\phi_j(x)=\sqrt{\frac{2j+1}{T-t}}P_j\biggl(\biggl(
x-\frac{T+t}{2}\biggr)\frac{2}{T-t}\biggr),\ \ \ j=0, 1, 2,\ldots,
\end{equation}

\vspace{5mm}
\noindent
where $P_j(x)$ is the Legendre polynomial. 

Using the system of 
functions (\ref{66}) and 
Theorems 1, 2 we obtain the following approximations of 
iterated 
It\^{o} stochastic integrals \cite{9}-\cite{1003} 
(also see early publications \cite{35} (1997), \cite{36} (1998),
\cite{37} (2000))

\vspace{1mm}
\begin{equation}
\label{opp0}
I_{(1)T,t}^{(i_1)}=\sqrt{T-t}\zeta_0^{(i_1)},
\end{equation}

\vspace{2mm}
\begin{equation}
\label{opp1}
I_{(01)T,t}^{(0 i_1)}=\frac{(T-t)^{3/2}}{2}\biggl(\zeta_0^{(i_1)}+
\frac{1}{\sqrt{3}}\zeta_1^{(i_1)}\biggr),\
\end{equation}

\vspace{2mm}
\begin{equation}
\label{opp2}
I_{(10)T,t}^{(i_1 0)}=\frac{(T-t)^{3/2}}{2}\biggl(\zeta_0^{(i_1)}-
\frac{1}{\sqrt{3}}\zeta_1^{(i_1)}\biggr),
\end{equation}

\vspace{4mm}
\begin{equation}
\label{kr00}
I_{(11)T,t}^{(i_1 i_2)q}=
\frac{T-t}{2}\left(\zeta_0^{(i_1)}\zeta_0^{(i_2)}+\sum_{i=1}^{q}
\frac{1}{\sqrt{4i^2-1}}\biggl(
\zeta_{i-1}^{(i_1)}\zeta_{i}^{(i_2)}-
\zeta_i^{(i_1)}\zeta_{i-1}^{(i_2)}\biggr)-{\bf 1}_{\{i_1=i_2\}}
\right),
\end{equation}

\vspace{7mm}
$$
I_{(111)T,t}^{(i_1 i_2 i_3)q_1}=
\sum_{j_1,j_2,j_3=0}^{q_1}
C_{j_3j_2j_1}
\Biggl(
\zeta_{j_1}^{(i_1)}\zeta_{j_2}^{(i_2)}\zeta_{j_3}^{(i_3)}
\Biggr.-{\bf 1}_{\{i_1=i_2\}}
{\bf 1}_{\{j_1=j_2\}}
\zeta_{j_3}^{(i_3)}-
$$

\vspace{2mm}
\begin{equation}
\label{kr1}
\Biggl.-{\bf 1}_{\{i_2=i_3\}}
{\bf 1}_{\{j_2=j_3\}}
\zeta_{j_1}^{(i_1)}-
{\bf 1}_{\{i_1=i_3\}}
{\bf 1}_{\{j_1=j_3\}}
\zeta_{j_2}^{(i_2)}\Biggr),
\end{equation}

\vspace{7mm}
$$
I_{(111)T,t}^{(i_1 i_1 i_1)}
=\frac{1}{6}(T-t)^{3/2}\biggl(
\biggl(\zeta_0^{(i_1)}\biggr)^3-3
\zeta_0^{(i_1)}\biggr),
$$

\vspace{7mm}

$$
C_{j_3j_2j_1}=\frac{\sqrt{(2j_1+1)(2j_2+1)(2j_3+1)}(T-t)^{3/2}}{8}\bar
C_{j_3j_2j_1},
$$

\vspace{4mm}
$$
\bar C_{j_3j_2j_1}=\int\limits_{-1}^{1}P_{j_3}(z)
\int\limits_{-1}^{z}P_{j_2}(y)
\int\limits_{-1}^{y}
P_{j_1}(x)dx dy dz,
$$

\vspace{6mm}
\noindent
standard Gaussian random variables 
$\zeta_{j}^{(i)}$ ($i\ne 0$) are defined by 
(\ref{rr23}), and

\vspace{2mm}
$$
I_{(11)T,t}^{(i_1 i_2)}=\hbox{\vtop{\offinterlineskip\halign{
\hfil#\hfil\cr
{\rm l.i.m.}\cr
$\stackrel{}{{}_{q\to \infty}}$\cr
}} }I_{(11)T,t}^{(i_1 i_2)q},
$$

\vspace{1mm}

$$
I_{(111)T,t}^{(i_1 i_2 i_3)}=\hbox{\vtop{\offinterlineskip\halign{
\hfil#\hfil\cr
{\rm l.i.m.}\cr
$\stackrel{}{{}_{q_1\to \infty}}$\cr
}} }I_{(111)T,t}^{(i_1 i_2 i_3)q_1}.
$$

\vspace{5mm}

Note that $T-t\ll 1$ ($T-t$ is an 
integration step with respect to the temporal variable). Thus
$q_1\ll q$ (see Table 1 \cite{9}-\cite{24aa},
\cite{arxiv-4}). 
Moreover, the values $\bar C_{j_3j_2j_1}$ 
do not depend on
$T-t.$ This feature is important because 
we can use a variable integration step $T-t$.
Coefficients 
$\bar C_{j_3j_2j_1}$ 
are calculated once and before the start 
of the numerical scheme.
Some examples of exact calculation of coefficients
$\bar C_{j_3j_2j_1}$ 
via Python programming language can be found in Table 2
(the database with 270,000 exactly
calculated Fourier--Legendre coefficients was described in \cite{1000}, \cite{1001}).

According to the notations introduced above, we have

$$
E^{(i_1i_2)q}={\sf M}\biggl\{\biggl(
I_{(11)T,t}^{(i_1 i_2)}-
I_{(11)T,t}^{(i_1 i_2)q}\biggr)^2\biggr\},
$$

\vspace{1mm}
$$
E^{(i_1i_2i_3)q_1}={\sf M}\biggl\{\biggl(
I_{(111)T,t}^{(i_1 i_2 i_3)}-
I_{(111)T,t}^{(i_1 i_2 i_3)q_1}\biggr)^2\biggr\}.
$$ 

\vspace{3mm}

\noindent
\begin{figure}
\begin{center}
\centerline{Table 1.\ Minimal numbers\ $q,$\ $q_1$\  
such that\ $E^{(i_1i_2)q},$\ $E^{(i_1i_2i_3)q_1}\le (T-t)^4,$\ \ \ $q_1\ll q$.}
\vspace{6mm}
\tabcolsep=0.28em
\begin{tabular}{|c|c|c|c|c|c|c|}
\hline
$T-t$&$0.08222$&$0.05020$&$0.02310$&$0.01956$\\
\hline
$q$&19&51&235&328\\
\hline
$q_1$&1&2&5&6\\
\hline
\end{tabular}
\end{center}
\vspace{5mm}
\begin{center}
\centerline{Table 2.\ Coefficients $\bar C_{3jk}.$}
\vspace{6mm}
\begin{tabular}{|c|c|c|c|c|c|c|c|c|}
\hline
${}_j {}^k$&0&1&2&3&4&5&6\\
\hline
0&$0$&$\frac{2}{105}$&$0$&$-\frac{4}{315}$&$0$&$\frac{2}{693}$&0\\
\hline
1&$\frac{4}{105}$&0&$-\frac{2}{315}$&0&$-\frac{8}{3465}$&0&$\frac{10}{9009}$\\
\hline
2&$\frac{2}{35}$&$-\frac{2}{105}$&$0$&$\frac{4}{3465}$&
$0$&$-\frac{74}{45045}$&0\\
\hline
3&$\frac{2}{315}$&$0$&$-\frac{2}{3465}$&0&
$\frac{16}{45045}$&0&$-\frac{10}{9009}$\\
\hline
4&$-\frac{2}{63}$&$\frac{46}{3465}$&0&$-\frac{32}{45045}$&
0&$\frac{2}{9009}$&0\\
\hline
5&$-\frac{10}{693}$&0&$\frac{38}{9009}$&0&
$-\frac{4}{9009}$&0&$\frac{122}{765765}$\\
\hline
6&$0$&$-\frac{10}{3003}$&$0$&$\frac{20}{9009}$&$0$&$-\frac{226}{765765}$&$0$\\
\hline
\end{tabular}
\end{center}
\vspace{5mm}
\end{figure}

\vspace{2mm}

Then for pairwise different $i_1,i_2,i_3=1,\ldots,m$
we obtain 
\cite{9}-\cite{24aa}, 
\cite{arxiv-3}

\vspace{1mm}
\begin{equation}
\label{909}
E^{(i_1i_2)q}=\frac{(T-t)^2}{2}\left(\frac{1}{2}-\sum_{i=1}^{q}
\frac{1}{4i^2-1}\right),
\end{equation}

\vspace{1mm}
\begin{equation}
\label{zzz}
E^{(i_1 i_2i_3)q_1}=
\frac{(T-t)^{3}}{6}-\sum_{j_1,j_2,j_3=0}^{q_1}
C_{j_3j_2j_1}^2.
\end{equation}

\vspace{4mm}

On the basis of 
the presented 
approximations of 
iterated It\^{o} stochastic integrals we 
can see that increasing of multiplicities of these integrals 
leads to increasing 
of orders of smallness with respect to $T-t$
($T-t\ll 1$) in the mean-square sense 
for iterated It\^{o} stochastic integrals. This leads to sharp decrease  
of member 
quantities
in the approximations of iterated It\^{o} stochastic 
integrals, which are required for achieving the acceptable accuracy
of approximation ($q_1\ll q$). 

From (\ref{zzz}) we obtain
\cite{9}-\cite{24aa}, \cite{arxiv-4}, \cite{arxiv-25}
(for more details see \cite{1000}-\cite{1003})

\vspace{2mm}
\begin{equation}
\label{ooo1}
\left.E^{(i_1i_2i_3)q_1}\right|_{q_1=6}
\approx
0.01956000(T-t)^3.
\end{equation}

\vspace{5mm}

\section{Approximation of Iterated Stochastic Integrals of Multiplicity
$k$ with
Respect to the $Q$-Wiener Process}

\vspace{5mm}

Consider the iterated It\^{o} stochastic integral with respect to 
the $Q$-Wiener process in the following form

\vspace{2mm}
$$
I[\Phi^{(k)}(Z), \psi^{(k)}]_{T,t}=
$$

\begin{equation}
\label{ssss11}
=
\int\limits_{t}^{T}\Phi_k(Z) \left( \ldots \left(
\int\limits_{t}^{t_3}\Phi_2(Z) \left(
\int\limits_{t}^{t_2}\Phi_1(Z)
\psi_1(t_1)d{\bf W}_{t_1} \right)
\psi_2(t_2)d{\bf W}_{t_2} \right) \ldots  \right) \psi_k(t_k)d{\bf W}_{t_k},\
\end{equation}

\vspace{8mm}
\noindent
where $Z: \Omega \rightarrow H$ is 
an ${\bf F}_t/{\mathcal{B}}(H)$-measurable mapping, 
$\Phi_k(v)(\ \ldots (\Phi_2(v)(\Phi_1(v))) \ldots\ )$ 
is a $k$-linear Hilbert--Schmidt operator
mapping from
$\underbrace{U_0\times \ldots \times U_0}_{\small{\hbox{$k$ times}}}$ to $H$
for all $v\in H,$ and $\psi_1(\tau),\ldots,\psi_k(\tau)\in L_2([t, T]).$

\vspace{1.5mm}

Let $I[\Phi^{(k)}(Z), \psi^{(k)}]_{T,t}^M$ be an approximation
of the stochastic integral (\ref{ssss11})

\vspace{3mm}
$$
I[\Phi^{(k)}(Z), \psi^{(k)}]_{T,t}^M=
$$

$$
=
\int\limits_{t}^{T}\Phi_k(Z) \left( \ldots \left(
\int\limits_{t}^{t_3}\Phi_2(Z) \left(
\int\limits_{t}^{t_2}\Phi_1(Z)
\psi_1(t_1)d{\bf W}_{t_1}^M \right)
\psi_2(t_2)d{\bf W}_{t_2}^M 
\right) \ldots  \right) \psi_k(t_k)d{\bf W}_{t_k}^M=
$$

\vspace{3mm}
$$
=\sum_{r_1,r_2,\ldots,r_k\in J_M}
\Phi_k(Z)\left(\ldots \left(\Phi_2(Z) \left(\Phi_1(Z)
e_{r_1} \right) 
e_{r_2} \right) \ldots \right) e_{r_k}\times
$$

\vspace{2mm}

\begin{equation}
\label{xx605}
\times
\sqrt{\lambda_{r_1}\lambda_{r_2}\ldots \lambda_{r_k}}\ 
J[\psi^{(k)}]_{T,t}^{(r_1 r_2\ldots r_k)},
\end{equation}

\vspace{7mm}
\noindent
where $0\le t<T\le \bar{T}$  and

\vspace{2mm}
$$
J[\psi^{(k)}]_{T,t}^{(r_1\ldots r_k)}=
\int\limits_t^T \psi_k(t_k)\ldots \int\limits_t^{t_3}
\psi_2(t_2)\int\limits_t^{t_{2}}\psi_1(t_1)
d{\bf w}_{t_1}^{(r_1)}d{\bf w}_{t_2}^{(r_2)}\ldots
d{\bf w}_{t_k}^{(r_k)}
$$

\vspace{5mm}
\noindent
is the iterated It\^{o} stochastic integral (\ref{sodom20}),
$r_1,r_2,\ldots,r_k\in J_M.$

Let $I[\Phi^{(k)}(Z), \psi^{(k)}]_{T,t}^{M,p_1\ldots,p_k}$ be an
approximation of the stochastic integral (\ref{xx605})

\vspace{5mm}
$$
I[\Phi^{(k)}(Z), \psi^{(k)}]_{T,t}^{M,p_1\ldots,p_k}
=\sum_{r_1,r_2,\ldots,r_k\in J_M}
\Phi_k(Z)\left(\ldots \left(\Phi_2(Z) \left(\Phi_1(Z)
e_{r_1} \right) 
e_{r_2} \right) \ldots \right) e_{r_k}\times
$$

\vspace{4mm}
\begin{equation}
\label{xx705}
\times
\sqrt{\lambda_{r_1}\lambda_{r_2}\ldots \lambda_{r_k}}\
J[\psi^{(k)}]_{T,t}^{(r_1 r_2\ldots r_k)p_1,\ldots,p_k},
\end{equation}

\vspace{6mm}
\noindent
where
$J[\psi^{(k)}]_{T,t}^{(r_1 r_2\ldots r_k)p_1,\ldots,p_k}$
is defined as a prelimit expression on the right-hand side of (\ref{leto6000})

\vspace{3mm}
$$
J[\psi^{(k)}]_{T,t}^{(r_1\ldots r_k)p_1\ldots p_k}=
\sum\limits_{j_1=0}^{p_1}\ldots
\sum\limits_{j_k=0}^{p_k}
C_{j_k\ldots j_1}\Biggl(
\prod_{l=1}^k\zeta_{j_l}^{(r_l)}+\sum\limits_{m=1}^{[k/2]}
(-1)^m \times
\Biggr.
$$

\vspace{2mm}
\begin{equation}
\label{f1}
\times
\sum_{\stackrel{(\{\{g_1, g_2\}, \ldots, 
\{g_{2m-1}, g_{2m}\}\}, \{q_1, \ldots, q_{k-2m}\})}
{{}_{\{g_1, g_2, \ldots, 
g_{2m-1}, g_{2m}, q_1, \ldots, q_{k-2m}\}=\{1, 2, \ldots, k\}}}}
\prod\limits_{s=1}^m
{\bf 1}_{\{r_{g_{{}_{2s-1}}}=~r_{g_{{}_{2s}}}\ne 0\}}
\Biggl.{\bf 1}_{\{j_{g_{{}_{2s-1}}}=~j_{g_{{}_{2s}}}\}}
\prod_{l=1}^{k-2m}\zeta_{j_{q_l}}^{(r_{q_l})}\Biggr).
\end{equation}

\vspace{7mm}

Let $U,$ $H$ be separable
$\mathbb{R}$-Hilbert spaces,
$U_{0}=Q^{1/2}(U)$, and
$L(U, H)$ be the space of linear and bounded operators mapping
from $U$ to $H$. Let
$L(U, H)_{0}=\left\{T \vert_{U_{0}}:\
T\in L(U, H)\right\}$ (here $T \vert_{U_{0}}$
is the restriction of operator $T$ to
the space $U_0$). It is known \cite{7} that
$L(U, H)_{0}$
is a dense subset of the space of
Hilbert--Schmidt operators $L_{HS}(U_{0}, H)$.

\vspace{2mm}

{\bf Theorem 3}\ \cite{9}, \cite{arxiv-20}, \cite{24a}-\cite{24aaaxxx}. 
{\it Suppose that
$\{\phi_j(x)\}_{j=0}^{\infty}$ is an arbitrary complete orthonormal system  
of functions in the space $L_2([t,T])$ and
$\psi_1(\tau),\ldots,\psi_k(\tau)\in L_2([t, T]).$ 
Furthermore$,$ let the following conditions be satisfied$:$

{\rm 1}. $Q\in L(U)$ is a nonnegative and symmetric 
trace class operator {\rm (}$\lambda_i$ and $e_i$ $(i\in J)$ are
its eigenvalues and eigenfunctions {\rm (}which form
an orthonormal basis of $U${\rm )}
correspondingly{\rm )}, and ${\bf W}_{\tau},$ $\tau\in [0, \bar T]$
is an $U$-valued $Q$-Wiener process.

{\rm 2}. $Z: \Omega \rightarrow H$ is 
an ${\bf F}_t/{\mathcal{B}}(H)$-measurable mapping.

{\rm 3}. $\Phi_1\in L(U, H)_{0},$ 
$\Phi_2\in L(H,L(U,H)_0),$ and $\Phi_k(v)(\ \ldots (\Phi_2(v)(\Phi_1(v))) \ldots\ )$ 
is a $k$-linear Hilbert--Schmidt operator mapping from
$\underbrace{U_0\times \ldots \times U_0}_{\small{\hbox{$k$ times}}}$ to $H$
for all $v\in H$
such that

\vspace{4mm}

$$
\Biggl\Vert \Phi_k(Z)\left(\ldots \left(\Phi_2(Z) \left(\Phi_1(Z)
e_{r_1} \right) 
e_{r_2} \right) \ldots \right) e_{r_k}\Biggr\Vert_H^2
\le L_k<\infty
$$

\vspace{6mm}
\noindent
w.\ p.\ {\rm 1}\ for all $r_1,r_2,\ldots,r_k\in J_M$, $M\in\mathbb{N}$.
Then

\vspace{2mm}
$$
{\sf M}\left\{
\Biggl\Vert
I[\Phi^{(k)}(Z), \psi^{(k)}]_{T,t}^M-
I[\Phi^{(k)}(Z), \psi^{(k)}]_{T,t}^{M,p_1\ldots p_k}
\Biggr\Vert_H^2\right\}\le 
$$

\vspace{1mm}
\begin{equation}
\label{zzz1}
\le L_k (k!)^2
\left({\rm tr}\ Q\right)^k
\left(I_k-\sum_{j_1=0}^{p_1}\ldots
\sum_{j_k=0}^{p_k}C^2_{j_k\ldots j_1}\right),
\end{equation}

\vspace{7mm}
\noindent
where $I_k$ is defined by {\rm (\ref{g123})}, 
${\rm tr}\ Q=\sum\limits_{i\in J}\lambda_i,$
and 

\vspace{1mm}
$$
C_{j_k\ldots j_1}=\int\limits_{[t,T]^k}
K(t_1,\ldots,t_k)\prod_{l=1}^{k}\phi_{j_l}(t_l)dt_1\ldots dt_k,
$$

\vspace{3mm}
$$
K(t_1,\ldots,t_k)=
\begin{cases}
\psi_1(t_1)\ldots \psi_k(t_k),\ t_1<\ldots<t_k\\
~\\
~\\
0,\ \hbox{\rm otherwise}
\end{cases}.
$$
}

\vspace{3mm}

{\bf Remark 2.}\ {\it It should be noted that the right-hand side 
of the inequality {\rm (\ref{zzz1})} is independent of $M$
and tends to zero if $p_1,\ldots,p_k\to\infty$ due to the
Parseval equality.}

\vspace{5mm}

\vspace{2mm}

\section{Approximation of Iterated Stochastic Integrals 
From the Exponential Milstein and Wagner--Platen Schemes for SPDEs}

\vspace{5mm}

This section is devoted to the approximation of iterated
stochastic integrals 
from the Milstein scheme (\ref{fff100})
and Wagner--Platen scheme (\ref{fff}) for SPDEs. These integrals 
have the
following form

\vspace{1mm}

\begin{equation}
\label{ret1}
J_1[B(Z)]_{T,t}=\int\limits_{t}^{T}B(Z)d{\bf W}_s,
\end{equation}

\begin{equation}
\label{ret2}
J_2[B(Z)]_{T,t}=A\left(\int\limits_{t}^{T}
\int\limits_{t}^{s}
B(Z)d{\bf W}_{\tau}ds-\frac{(T-t)}{2}
\int\limits_{t}^{T}B(Z)d{\bf W}_s\right),
\end{equation}

\vspace{2mm}
\begin{equation}
\label{ret3}
J_3[B(Z),F(Z)]_{T,t}=
\int\limits_{t}^{T}B'(Z) \left(
\int\limits_{t}^{t_2}\biggl(AZ+F(Z)\biggr)dt_1 \right) d{\bf W}_{t_2},
\end{equation}

\vspace{2mm}
\begin{equation}
\label{ret4}
J_4[B(Z),F(Z)]_{T,t}=\int\limits_{t}^{T}F'(Z)
\left(\int\limits_{t}^{t_2} B(Z) d{\bf W}_{t_1} \right) dt_2,
\end{equation}

\vspace{2mm}
\begin{equation}
\label{ccc1}
I_1[B(Z)]_{T,t}=\int\limits_{t}^{T}
B'(Z)\left(\int\limits_{t}^s
B(Z)d{\bf W}_{\tau}\right)d{\bf W}_s,
\end{equation}

\vspace{2mm}
\begin{equation}
\label{ccc2}
I_2[B(Z)]_{T,t}=\int\limits_{t}^{T}B'(Z) \left(
\int\limits_{t}^{t_3}B'(Z)\left(
\int\limits_{t}^{t_2}B(Z) d{\bf W}_{t_1} 
\right) d{\bf W}_{t_2} \right) d{\bf W}_{t_3},
\end{equation}

\vspace{2mm}
\begin{equation}
\label{ccc3}
I_3[B(Z)]_{T,t}=\int\limits_{t}^{T}B''(Z) \left(
\int\limits_{t}^{t_2}B(Z) d{\bf W}_{t_1},
\int\limits_{t}^{t_2}B(Z) d{\bf W}_{t_1} 
\right) d{\bf W}_{t_2},
\end{equation}

\vspace{5mm}
\noindent
where $Z: \Omega \rightarrow H$ is 
an ${\bf F}_t/{\mathcal{B}}(H)$-measurable mapping,
$0\le t<T\le \bar T.$

Note that according to
(\ref{ggg0})--(\ref{ggg3}), (\ref{opp0}), and (\ref{opp1})
we can write w. p. 1

\vspace{2mm}
$$
J_1[B(Z)]_{T,t}^M=\int\limits_{t}^{T}B(Z)d{\bf W}_s^M=
(T-t)^{1/2}\sum\limits_{r_1\in J_M}
B(Z)e_{r_1}\sqrt{\lambda_{r_1}}
\zeta_0^{(r_1)},
$$

\vspace{7mm}

$$
J_2[B(Z)]_{T,t}^M=A\left(\int\limits_{t}^{T}
\int\limits_{t}^{s}
B(Z)d{\bf W}_{\tau}^M ds-\frac{(T-t)}{2}
\int\limits_{t}^{T}B(Z)d{\bf W}_s^M\right)=
$$

\vspace{1mm}
\begin{equation}
\label{sd11}
=-\frac{(T-t)^{3/2}}{2\sqrt{3}}\sum\limits_{r_1\in J_M}
AB(Z)e_{r_1}\sqrt{\lambda_{r_1}}
\zeta_1^{(r_1)},
\end{equation}

\vspace{7mm}

$$
J_3[B(Z),F(Z)]_{T,t}^M=
\int\limits_{t}^{T}B'(Z) \left(
\int\limits_{t}^{t_2}\biggl(AZ+F(Z)\biggr)dt_1 \right) d{\bf W}_{t_2}^M=
$$

\vspace{2mm}
\begin{equation}
\label{sd12}
=
\frac{(T-t)^{3/2}}{2}\sum\limits_{r_1\in J_M}
B'(Z)\biggl(AZ+F(Z)\biggr)
e_{r_1}\sqrt{\lambda_{r_1}}\left(\zeta_0^{(r_1)}+
\frac{1}{\sqrt{3}}\zeta_1^{(r_1)}\right),
\end{equation}

\vspace{7mm}

$$
J_4[B(Z),F(Z)]_{T,t}^M=\int\limits_{t}^{T}F'(Z)
\left(\int\limits_{t}^{t_2} B(Z) d{\bf W}_{t_1}^M \right) dt_2=
$$

\vspace{2mm}
\begin{equation}
\label{sd13}
=\frac{(T-t)^{3/2}}{2}\sum\limits_{r_1\in J_M}
F'(Z)B(Z)e_{r_1}\sqrt{\lambda_{r_1}}
\left(\zeta_0^{(r_1)}-
\frac{1}{\sqrt{3}}\zeta_1^{(r_1)}\right),
\end{equation}

\vspace{5mm}
\noindent
where $\zeta_0^{(r_1)},$ $\zeta_1^{(r_1)}$ $(r_1\in J_M)$ are 
independent standard Gaussian random variables.

Let $I_1[B(Z)]_{T,t}^M,$ $I_2[B(Z)]_{T,t}^M,$ 
$I_3[B(Z)]_{T,t}^M$ be approximations
of stochastic integrals (\ref{ccc1})--(\ref{ccc3}),
which have the following form (see (\ref{ch100}), (\ref{ch101}), and (\ref{da1}))

\vspace{2mm}
$$
I_1[B(Z)]_{T,t}^M=\int\limits_{t}^{T}
B'(Z)\left(\int\limits_{t}^s
B(Z)d{\bf W}_{\tau}^M\right)d{\bf W}_s^M=
$$

\vspace{2mm}
\begin{equation}
\label{vvv1}
=\sum\limits_{r_1,r_2\in J_M}
B'(Z)\left(B(Z)e_{r_1}\right)
e_{r_2}\sqrt{\lambda_{r_1}\lambda_{r_2}}I_{(11)T,t}^{(r_1 r_2)},
\end{equation}

\vspace{7mm}

$$
I_2[B(Z)]_{T,t}^M=
\int\limits_{t}^{T}B'(Z) \left(
\int\limits_{t}^{t_3}B'(Z)\left(
\int\limits_{t}^{t_2}B(Z) d{\bf W}_{t_1}^M 
\right) d{\bf W}_{t_2}^M \right) d{\bf W}_{t_3}^M=
$$

\vspace{2mm}
\begin{equation}
\label{da4xx}
=\sum_{r_1,r_2,r_3\in J_M}B'(Z)\left(B'(Z)\left(B(Z)e_{r_1}\right)
e_{r_2}\right)e_{r_3}
\sqrt{\lambda_{r_1}\lambda_{r_2}\lambda_{r_3}}
I_{(111)T,t}^{(r_1 r_2 r_3)},
\end{equation}

\vspace{7mm}

$$
I_3[B(Z)]_{T,t}^M=\int\limits_{t}^{T}B''(Z) \left(
\int\limits_{t}^{t_2}B(Z) d{\bf W}_{t_1}^M,
\int\limits_{t}^{t_2}B(Z) d{\bf W}_{t_1}^M 
\right) d{\bf W}_{t_2}^M=
$$

\vspace{2mm}
$$
=\sum_{r_1,r_2,r_3\in J_M}B''(Z)\left(B(Z)e_{r_1}, B(Z)e_{r_2}\right)e_{r_3}
\sqrt{\lambda_{r_1}\lambda_{r_2}\lambda_{r_3}} \times
$$

\vspace{2mm}
\begin{equation}
\label{da4x}
\times
\left(I_{(111)T,t}^{(r_1 r_2 r_3)}+
I_{(111)T,t}^{(r_2 r_1 r_3)}+
{\bf 1}_{\{r_1=r_2\}}I_{(01)T,t}^{(0 r_3)}\right).
\end{equation}

\vspace{9mm}

Let $I_1[B(Z)]_{T,t}^{M,q},$
$I_2[B(Z)]_{T,t}^{M.q},$ $I_3[B(Z)]_{T,t}^{M,q}$ be approximations
of stochastic integrals (\ref{vvv1})--(\ref{da4x}),
which are represented as follows

\vspace{4mm}

$$
I_1[B(Z)]_{T,t}^{M,q}=\sum\limits_{r_1,r_2\in J_M}
B'(Z)\left(B(Z)e_{r_1}\right)
e_{r_2}\sqrt{\lambda_{r_1}\lambda_{r_2}}I_{(11)T,t}^{(r_1 r_2)q},
$$

\vspace{8mm}

$$
I_2[B(Z)]_{T,t}^{M,q}=
\sum_{r_1,r_2,r_3\in J_M}B'(Z)\left(B'(Z)\left(B(Z)e_{r_1}\right)
e_{r_2}\right)e_{r_3}\times
$$

\vspace{2mm}
\begin{equation}
\label{da4xxx}
\times
\sqrt{\lambda_{r_1}\lambda_{r_2}\lambda_{r_3}}
I_{(111)T,t}^{(r_1 r_2 r_3)q},
\end{equation}

\vspace{6mm}

$$
I_3[B(Z)]_{T,t}^{M,q}=
\sum_{r_1,r_2,r_3\in J_M}B''(Z)\left(B(Z)e_{r_1}, B(Z)e_{r_2}\right)e_{r_3}
\sqrt{\lambda_{r_1}\lambda_{r_2}\lambda_{r_3}} \times
$$

\vspace{2mm}
\begin{equation}
\label{da4}
\times
\left(I_{(111)T,t}^{(r_1 r_2 r_3)q}+
I_{(111)T,t}^{(r_2 r_1 r_3)q}+
{\bf 1}_{\{r_1=r_2\}}I_{(01)T,t}^{(0 r_3)q}\right),
\end{equation}

\vspace{4mm}
\noindent
where $q\ge 1$ and the approximations 
$I_{(11)T,t}^{(r_1 r_2)q},$ $I_{(111)T,t}^{(r_1 r_2 r_3)q},$
$I_{(111)T,t}^{(r_2 r_1 r_3)q}$ are defined by 
(\ref{kr00}), (\ref{kr1}).

Recall that $L_{HS}(U_0,H)$
is a space of Hilbert--Schmidt operators mapping from $U_0$ to $H.$
Let $L_{HS}^{(2)}(U_0, H)$ and
$L_{HS}^{(3)}(U_0, H)$
be spaces of bilinear and 3-linear
Hilbert--Schmidt operators 
mapping 
from $U_0 \times U_0$ to $H$ and from $U_0 \times U_0 \times U_0$ to $H$
correspondingly.
Furthermore, let
$\left\Vert \cdot \right\Vert_{L_{HS}(U_0,H)}$,
$\left\Vert \cdot \right\Vert_{L_{HS}^{(2)}(U_0, H)},$
and
$\left\Vert \cdot \right\Vert_{L_{HS}^{(3)}(U_0, H)}$ 
be operator norms in these spaces.

\vspace{3mm}

{\bf Theorem 4}\ \cite{1004} (also see \cite{24a}-\cite{24aaaxxx}, \cite{1004a}). {\it 
Let the conditions {\rm 1, 2} of Theorem {\rm 3} be fulfilled.
Let $B(v)$ be a Hilbert--Schmidt operator mapping from
$U_0$ to $H$ for all $v\in H,$
$B'(v)(B(v))$ be
a bilinear 
Hilbert--Schmidt operator mapping from
$U_0 \times U_0$ to $H$
for all $v\in H,$ and
$B'(v)(B'(v)(B(v))),$ $B''(v)(B(v),B(v))$ be
{\rm 3}-linear 
Hilbert--Schmidt operators mapping from
$U_0\times U_0 \times U_0$ to $H$
for all $v\in H$ {\rm (}we suppose that Fr\^{e}chet derivatives 
$B',$ $B''$ exist {\rm (}see Sect.~{\rm 2))}. 
Moreover, 
let there exists a constant $C$ such that w. p. {\rm 1}

\vspace{2mm}
$$
\biggl\Vert B(Z)Q^{-\alpha}\biggr\Vert_{L_{HS}(U_0,H)}<C,\ \ \
\biggl\Vert B'(Z)(B(Z))Q^{-\alpha}\biggr\Vert_{L_{HS}^{(2)}(U_0,H)}<C
$$

\vspace{4mm}
$$
\biggl\Vert B'(Z)(B'(Z)(B(Z)))Q^{-\alpha}
\biggr\Vert_{L_{HS}^{(3)}(U_0,H)}<C,
$$

\vspace{4mm}
$$
\biggl\Vert B''(Z)(B(Z), B(Z))Q^{-\alpha}
\biggr\Vert_{L_{HS}^{(3)}(U_0,H)}<C,
$$

\vspace{7mm}
\noindent
for some $\alpha>0.$
Then

\vspace{2mm}
$$
{\sf M}\left\{\Biggl\Vert
I_1[B(Z)]_{T,t}-
I_1[B(Z)]_{T,t}^{M,q}\Biggr\Vert_H^2\right\}\le
$$

\begin{equation}
\label{rock6}
\le (T-t)^2\left(C_0
\left({\rm tr}\ Q\right)^2
\left(\frac{1}{2}-\sum_{j=1}^{q}
\frac{1}{4j^2-1}\right)+
K_Q\left(\sup\limits_{i\in J\backslash J_M}\lambda_i\right)^{2\alpha}
\right),
\end{equation}

\vspace{8mm}

$$
{\sf M}\left\{\Biggl\Vert
I_2[B(Z)]_{T,t}-
I_2[B(Z)]_{T,t}^{M,q}\Biggr\Vert_H^2\right\}\le
$$

\begin{equation}
\label{ppp6}
\le (T-t)^3\left(C_1
\left({\rm tr}\ Q\right)^3
\Biggl(\frac{1}{6}-\sum_{j_1,j_2,j_3=0}^{q}
{\hat C}_{j_3j_2j_1}^2\Biggr)+
L_Q\left(\sup\limits_{i\in J\backslash J_M}\lambda_i\right)^{2\alpha}
\right),
\end{equation}

\vspace{9mm}

$$
{\sf M}\left\{\Biggl\Vert
I_3[B(Z)]_{T,t}-
I_3[B(Z)]_{T,t}^{M,q}\Biggr\Vert_H^2\right\}\le
$$

\begin{equation}
\label{uuu1}
\le (T-t)^3\left(C_2
\left({\rm tr}\ Q\right)^3
\Biggl(\frac{1}{6}-\sum_{j_1,j_2,j_3=0}^{q}
{\hat C}_{j_3j_2j_1}^2\Biggr)+
M_Q\left(\sup\limits_{i\in J\backslash J_M}\lambda_i\right)^{2\alpha}
\right),
\end{equation}

\vspace{8mm}
\noindent
where $q\in\mathbb{N},$\ $C_0, C_1, C_2, K_Q, L_Q, M_Q<\infty,$ and 

\vspace{4mm}
$$
{\hat C}_{j_3j_2j_1}=\frac{\sqrt{(2j_1+1)(2j_2+1)(2j_3+1)}}{8}\bar
C_{j_3j_2j_1},
$$

\vspace{2mm}
$$
\bar C_{j_3j_2j_1}=\int\limits_{-1}^{1}P_{j_3}(z)
\int\limits_{-1}^{z}P_{j_2}(y)
\int\limits_{-1}^{y}
P_{j_1}(x)dx dy dz,
$$

\vspace{5mm}
\noindent 
where $P_j(x)\ (j=0, 1, 2,\ldots)$ is the Legendre polynomial.
}

\vspace{2mm}

{\bf Remark 3.}\ {\it Note that the estimate like 
{\rm (\ref{rock6})} has been derived in \cite{8} 
{\rm (}also see \cite{2}{\rm )}
with the 
difference connected with the first term
on the right-hand side of {\rm (\ref{rock6})}.
In \cite{8} the authors used the Karhunen--Loeve
expansion of the Brownian bridge process for the
approximation of iterated It\^{o} stochastic integrals
with respect to the scalar standard Wiener processes.
In this article we apply Theorem {\rm 1} and the
system of Legendre polynomials for obtainment
the first term on the right-hand side of {\rm (\ref{rock6}).}}

\vspace{2mm}

{\bf Proof.}\ The estimate (\ref{rock6}) directly follows
from Theorem 3 of this article
(the first term on the right-hand side of (\ref{rock6}))
and Theorem 1 from \cite{8} (the second term on the 
right-hand side of (\ref{rock6})).
Further $C_3,$ $C_4,$ $\ldots $ denote various constants.

Let us prove the estimates (\ref{ppp6}), (\ref{uuu1}).
Using the elementary inequality $(a+b)^2\le 2(a^2+b^2)$ and Theorem 3,
we obtain

\vspace{1mm}

$$
{\sf M}\left\{\Biggl\Vert
I_2[B(Z)]_{T,t}-
I_2[B(Z)]_{T,t}^{M,q}\Biggr\Vert_H^2\right\}\le 
$$

\vspace{2mm}
$$
\le 2\left(
{\sf M}\left\{\Biggl\Vert
I_2[B(Z)]_{T,t}-
I_2[B(Z)]_{T,t}^{M}\Biggr\Vert_H^2\right\}+
{\sf M}\left\{\Biggl\Vert
I_2[B(Z)]_{T,t}^M-
I_2[B(Z)]_{T,t}^{M,q}\Biggr\Vert_H^2\right\}\right)\le
$$

\vspace{2mm}
\begin{equation}
\label{ppp4}
\le 
2{\sf M}\left\{\Biggl\Vert
I_2[B(Z)]_{T,t}-
I_2[B(Z)]_{T,t}^{M}\Biggr\Vert_H^2\right\}+
C_3(T-t)^3
\left({\rm tr}\ Q\right)^3
\Biggl(\frac{1}{6}-\sum_{j_1,j_2,j_3=0}^{q}
{\hat C}_{j_3j_2j_1}^2\Biggr),
\end{equation}

\vspace{7mm}

$$
{\sf M}\left\{\Biggl\Vert
I_3[B(Z)]_{T,t}-
I_3[B(Z)]_{T,t}^{M,q}\Biggr\Vert_H^2\right\}\le 
$$

\vspace{2mm}
\begin{equation}
\label{hhh}
\le 2\left(
{\sf M}\left\{\Biggl\Vert
I_3[B(Z)]_{T,t}-
I_3[B(Z)]_{T,t}^{M}\Biggr\Vert_H^2\right\}+
{\sf M}\left\{\Biggl\Vert
I_3[B(Z)]_{T,t}^M-
I_3[B(Z)]_{T,t}^{M,q}\Biggr\Vert_H^2\right\}\right).
\end{equation}

\vspace{6mm}

Repeating with an insignificant 
modification the proof 
of Theorem 3 for the case $k=3$ (see for delails 
\cite{9} (pp. 39--44) or \cite{arxiv-20}, \cite{24a}-\cite{24aaaxxx}), 
we have

\vspace{1mm}
\begin{equation}
\label{hhh1}
{\sf M}\left\{\biggl\Vert
I_3[B(Z)]_{T,t}^{M}-I_3[B(Z)]_{T,t}^{M,q}\biggr\Vert_H^2\right\}\le
4C(3!)^2
\left({\rm tr}\ Q\right)^3 (T-t)^3
\Biggl(\frac{1}{6}-\sum_{j_1,j_2,j_3=0}^{q}
{\hat C}_{j_3j_2j_1}^2\Biggr),
\end{equation}

\vspace{4mm}
\noindent
where constant $C$ has the same meaning 
as constant $L_k$ in Theorem 3 
($k$ is the multiplicity of the iterated stochastic integral).

Combining (\ref{hhh}) and (\ref{hhh1}), we obtain

\vspace{2mm}
$$
{\sf M}\left\{\Biggl\Vert
I_3[B(Z)]_{T,t}-
I_3[B(Z)]_{T,t}^{M,q}\Biggr\Vert_H^2\right\}\le 
$$

\vspace{2mm}
\begin{equation}
\label{sss1}
\le 
2{\sf M}\left\{\Biggl\Vert
I_3[B(Z)]_{T,t}-
I_3[B(Z)]_{T,t}^{M}\Biggr\Vert_H^2\right\}+
C_4(T-t)^3
\left({\rm tr}\ Q\right)^3
\Biggl(\frac{1}{6}-\sum_{j_1,j_2,j_3=0}^{q}
{\hat C}_{j_3j_2j_1}^2\Biggr).
\end{equation}

\vspace{7mm}

Let us evaluate the values

\vspace{2mm}
$$
{\sf M}\left\{\Biggl\Vert
I_2[B(Z)]_{T,t}-
I_2[B(Z)]_{T,t}^{M}\Biggr\Vert_H^2\right\},\ \ \
{\sf M}\left\{\Biggl\Vert
I_3[B(Z)]_{T,t}-
I_3[B(Z)]_{T,t}^{M}\Biggr\Vert_H^2\right\}.
$$

\vspace{7mm}

Using the elementary inequality 
$(a+b+c)^2\le 3(a^2+b^2+c^2)$, we have

\vspace{2mm}

\begin{equation}
\label{ppp0}
{\sf M}\left\{\Biggl\Vert
I_2[B(Z)]_{T,t}-
I_2[B(Z)]_{T,t}^{M}\Biggr\Vert_H^2\right\}\le
3\left(E^{1,M}_{T,t}+E^{2,M}_{T,t}+E^{3,M}_{T,t}\right),
\end{equation}

\vspace{2mm}
\begin{equation}
\label{ppp0x}
{\sf M}\left\{\Biggl\Vert
I_3[B(Z)]_{T,t}-
I_3[B(Z)]_{T,t}^{M}\Biggr\Vert_H^2\right\}\le
3\left(G^{1,M}_{T,t}+G^{2,M}_{T,t}+G^{3,M}_{T,t}\right),
\end{equation}

\vspace{6mm}
\noindent 
where

\vspace{2mm}
$$
E^{1,M}_{T,t}=
{\sf M}\left\{\left\Vert
\int\limits_{t}^{T}B'(Z) \left(
\int\limits_{t}^{t_3}B'(Z)\left(
\int\limits_{t}^{t_2}B(Z) d\left({\bf W}_{t_1}-{\bf W}_{t_1}^M\right) 
\right) d{\bf W}_{t_2} \right) d{\bf W}_{t_3}\right\Vert_H^2\right\},
$$

\vspace{4mm}
$$
E^{2,M}_{T,t}=
{\sf M}\left\{\left\Vert
\int\limits_{t}^{T}B'(Z) \left(
\int\limits_{t}^{t_3}B'(Z)\left(
\int\limits_{t}^{t_2}B(Z) d{\bf W}_{t_1}^M
\right) d\left({\bf W}_{t_2}-{\bf W}_{t_2}^M\right) \right) 
d{\bf W}_{t_3}\right\Vert_H^2\right\},
$$

\vspace{4mm}
$$
E^{3,M}_{T,t}=
{\sf M}\left\{\left\Vert
\int\limits_{t}^{T}B'(Z) \left(
\int\limits_{t}^{t_3}B'(Z)\left(
\int\limits_{t}^{t_2}B(Z) d{\bf W}_{t_1}^M
\right) d{\bf W}_{t_2}^M\right) 
d\left({\bf W}_{t_3}-{\bf W}_{t_3}^M\right)\right\Vert_H^2\right\},
$$

\vspace{4mm}
$$
G^{1,M}_{T,t}=
{\sf M}\left\{\left\Vert
\int\limits_{t}^{T}B''(Z) \left(
\int\limits_{t}^{t_2}B(Z)d{\bf W}_{t_1},
\int\limits_{t}^{t_2}B(Z)d\left({\bf W}_{t_1}-{\bf W}_{t_1}^M\right)
\right) 
d{\bf W}_{t_2}\right\Vert_H^2\right\},
$$

\vspace{4mm}
$$
G^{2,M}_{T,t}=
{\sf M}\left\{\left\Vert
\int\limits_{t}^{T}B''(Z) \left(
\int\limits_{t}^{t_2}B(Z)d\left({\bf W}_{t_1}-{\bf W}_{t_1}^M\right),
\int\limits_{t}^{t_2}B(Z)d{\bf W}_{t_1}^M
\right) 
d{\bf W}_{t_2}\right\Vert_H^2\right\},
$$

\vspace{4mm}
$$
G^{3,M}_{T,t}=
{\sf M}\left\{\left\Vert
\int\limits_{t}^{T}B''(Z) \left(
\int\limits_{t}^{t_2}B(Z)d{\bf W}_{t_1}^M,
\int\limits_{t}^{t_2}B(Z)d{\bf W}_{t_1}^M
\right) 
d\left({\bf W}_{t_2}-{\bf W}_{t_2}^M\right)\right\Vert_H^2\right\}.
$$

\vspace{9mm}

We have

\vspace{2mm}
$$
E^{1,M}_{T,t}=
\int\limits_{t}^{T}{\sf M}\left\{\left\Vert
B'(Z) \left(
\int\limits_{t}^{t_3}B'(Z)\left(
\int\limits_{t}^{t_2}B(Z) d\left({\bf W}_{t_1}-{\bf W}_{t_1}^M\right) 
\right) d{\bf W}_{t_2} \right)\right\Vert_{L_{HS}(U_0,H)}^2
\right\}dt_3\le 
$$

\vspace{4mm}
$$
\le C_5 
\int\limits_{t}^{T}{\sf M}\left\{\left\Vert
\int\limits_{t}^{t_3}B'(Z)\left(
\int\limits_{t}^{t_2}B(Z) d\left({\bf W}_{t_1}-{\bf W}_{t_1}^M\right) 
\right) d{\bf W}_{t_2}\right\Vert_H^2
\right\}dt_3=
$$

\vspace{4mm}
$$
=
C_5 
\int\limits_{t}^{T}\int\limits_{t}^{t_3}{\sf M}\left\{\left\Vert
B'(Z)\left(
\int\limits_{t}^{t_2}B(Z) d\left({\bf W}_{t_1}-{\bf W}_{t_1}^M\right) 
\right)\right\Vert_{L_{HS}(U_0,H)}^2 \right\}dt_2
dt_3\le
$$

\vspace{4mm}
\begin{equation}
\label{22u}
\le 
C_6 
\int\limits_{t}^{T}\int\limits_{t}^{t_3}{\sf M}\left\{\left\Vert
\int\limits_{t}^{t_2}B(Z) d\left({\bf W}_{t_1}-{\bf W}_{t_1}^M\right) 
\right\Vert_H^2 \right\}dt_2
dt_3\le
\end{equation}

\vspace{4mm}
\begin{equation}
\label{23u}
\le 
C_6 \left(\sup\limits_{i\in J\backslash J_M}\lambda_i\right)^{2\alpha}
\int\limits_{t}^{T}\int\limits_{t}^{t_3}
\int\limits_{t}^{t_2}
{\sf M}\left\{\biggl\Vert
B(Z)Q^{-\alpha}\biggr\Vert_{L_{HS}(U_0,H)}^2 \right\}dt_1 dt_2
dt_3 \le 
\end{equation}

\vspace{3mm}
\begin{equation}
\label{ppp1}
\le C_7 
\left(\sup\limits_{i\in J\backslash J_M}\lambda_i\right)^{2\alpha}
(T-t)^3.
\end{equation}

\vspace{8mm}

Note that the transition from (\ref{22u}) to (\ref{23u}) was made
by analogy with the proof of Theorem 1 in \cite{8} (also see \cite{2}). 
More precisely, taking into account the relation $Q^{\alpha}e_i=
\lambda_i^{\alpha}e_i$, we have 
(see \cite{8}, Sect. 3.1)

$$
{\sf M}\left\{\left\Vert
\int\limits_{t}^{t_2}B(Z) d\left({\bf W}_{t_1}-{\bf W}_{t_1}^M\right) 
\right\Vert_H^2\right\}=
$$

\vspace{4mm}
$$
=
{\sf M}\left\{\left\Vert \sum\limits_{i\in J\backslash J_M}
\sqrt{\lambda_i}
\int\limits_{t}^{t_2}B(Z)e_i d{\bf w}_{t_1}^{(i)} 
\right\Vert_H^2\right\}=
$$

\vspace{4mm}
$$
=\sum\limits_{i\in J\backslash J_M}\lambda_i
\int\limits_{t}^{t_2}{\sf M}\left\{\biggl\Vert
B(Z)Q^{-\alpha}Q^{\alpha}e_i 
\biggr\Vert_H^2\right\}dt_1=
$$

\vspace{4mm}
$$
=\sum\limits_{i\in J\backslash J_M}\lambda_i^{1+2\alpha}
\int\limits_{t}^{t_2}{\sf M}\left\{\biggl\Vert
B(Z)Q^{-\alpha}e_i 
\biggr\Vert_H^2\right\}dt_1=
$$

\vspace{4mm}
$$
=\left(\sup\limits_{i\in J\backslash J_M}\lambda_i\right)^{2\alpha}
\int\limits_{t}^{t_2}{\sf M}\left\{
\sum\limits_{i\in J\backslash J_M}\lambda_i
\biggl\Vert
B(Z)Q^{-\alpha}e_i 
\biggr\Vert_H^2\right\}dt_1\le
$$

\vspace{4mm}
$$
\le 
\left(\sup\limits_{i\in J\backslash J_M}\lambda_i\right)^{2\alpha}
\int\limits_{t}^{t_2}{\sf M}\left\{
\sum\limits_{i\in J}\lambda_i
\biggl\Vert
B(Z)Q^{-\alpha}e_i 
\biggr\Vert_H^2\right\}dt_1=
$$

\vspace{4mm}
\begin{equation}
\label{24u}
=\left(\sup\limits_{i\in J\backslash J_M}\lambda_i\right)^{2\alpha}
\int\limits_{t}^{t_2}{\sf M}\left\{
\biggl\Vert
B(Z)Q^{-\alpha} 
\biggr\Vert_{L_{HS}(U_0,H)}^2\right\}dt_1.
\end{equation}

\vspace{8mm}

Further
we also will use the estimate like (\ref{24u}).
We have

\vspace{3mm}
$$
E^{2,M}_{T,t}=
$$

\vspace{1mm}
$$
=
\int\limits_{t}^{T}{\sf M}\left\{\left\Vert
B'(Z) \left(
\int\limits_{t}^{t_3}B'(Z)\left(
\int\limits_{t}^{t_2}B(Z) d{\bf W}_{t_1}^M
\right) d\left({\bf W}_{t_2}-{\bf W}_{t_2}^M\right) 
\right)\right\Vert_{L_{HS}(U_0,H)}^2
\right\}dt_3\le 
$$

\vspace{4mm}
$$
\le C_8
\int\limits_{t}^{T}{\sf M}\left\{\left\Vert
\int\limits_{t}^{t_3}B'(Z)\left(
\int\limits_{t}^{t_2}B(Z) d{\bf W}_{t_1}^M
\right) d\left({\bf W}_{t_2}-{\bf W}_{t_2}^M\right) 
\right\Vert_H^2
\right\}dt_3\le 
$$

\vspace{4mm}
$$
\le C_8 \left(\sup\limits_{i\in J\backslash J_M}\lambda_i\right)^{2\alpha}
\int\limits_{t}^{T}\int\limits_{t}^{t_3}{\sf M}\left\{\left\Vert
B'(Z)\left(
\int\limits_{t}^{t_2}B(Z) d{\bf W}_{t_1}^M
\right)Q^{-\alpha}
\right\Vert_{L_{HS}(U_0,H)}^2
\right\}dt_2 dt_3 \le 
$$

\vspace{4mm}
$$
\le C_{8} \left(\sup\limits_{i\in J\backslash J_M}\lambda_i\right)^{2\alpha}
\int\limits_{t}^{T}\int\limits_{t}^{t_3}
{\sf M}\left\{\biggl\Vert B'(Z)\left(B(Z)\right)Q^{-\alpha}
\biggr\Vert^2_{L_{HS}^{(2)}(U_0,H)}\right\}(t_2-t)dt_2
dt_3 \le
$$

\vspace{3mm}
\begin{equation}
\label{ppp2}
\le C_9 
\left(\sup\limits_{i\in J\backslash J_M}\lambda_i\right)^{2\alpha}
(T-t)^3,
\end{equation}

\vspace{11mm}

$$
E^{3,M}_{T,t}\le
\left(\sup\limits_{i\in J\backslash J_M}\lambda_i\right)^{2\alpha}\times
$$

\vspace{4mm}
$$
\times
\int\limits_{t}^{T}{\sf M}\left\{\left\Vert
B'(Z) \left(
\int\limits_{t}^{t_3}B'(Z)\left(
\int\limits_{t}^{t_2}B(Z) d{\bf W}_{t_1}^M
\right) d{\bf W}_{t_2}^M\right) 
Q^{-\alpha}\right\Vert_{L_{HS}(U_0,H)}^2\right\}dt_3\le
$$

\vspace{4mm}
$$
\le C_{10}
\left(\sup\limits_{i\in J\backslash J_M}\lambda_i\right)^{2\alpha}
\int\limits_{t}^{T}
{\sf M}\left\{\biggl\Vert B'(Z)\left(B'(Z)\left(B(Z)\right)
\right)Q^{-\alpha}
\biggr\Vert^2_{L_{HS}^{(3)}(U_0,H)}\right\}
\frac{(t_3-t)^2}{2}dt_3
\le
$$

\vspace{2mm}
\begin{equation}
\label{ppp3}
\le C_{11}\left(\sup\limits_{i\in J\backslash J_M}\lambda_i\right)^{2\alpha}
(T-t)^3.
\end{equation}

\vspace{10mm}

Combining (\ref{ppp4}), (\ref{ppp0}), (\ref{ppp1})--(\ref{ppp3}),
we obtain (\ref{ppp6}).
We have

\vspace{1mm}
$$
G^{1,M}_{T,t}=\int\limits_{t}^{T}
{\sf M}\left\{\left\Vert
B''(Z) \left(
\int\limits_{t}^{t_2}B(Z)d{\bf W}_{t_1},
\int\limits_{t}^{t_2}B(Z)d\left({\bf W}_{t_1}-{\bf W}_{t_1}^M\right)
\right)\right\Vert_{L_{HS}(U_0,H)}^2\right\}dt_3\le
$$

\vspace{4mm}
$$
\le C_{12}\int\limits_{t}^{T}
{\sf M}\left\{\left\Vert
\int\limits_{t}^{t_2}B(Z)d{\bf W}_{t_1}
\right\Vert_H^2
\left\Vert
\int\limits_{t}^{t_2}B(Z)d\left({\bf W}_{t_1}-{\bf W}_{t_1}^M\right)
\right\Vert_H^2\right\}dt_3\le
$$

\vspace{4mm}
$$
\le C_{12}\int\limits_{t}^{T}
\left({\sf M}\left\{\left\Vert
\int\limits_{t}^{t_2}B(Z)d{\bf W}_{t_1}
\right\Vert_H^4\right\}\right)^{1/2}
\left({\sf M}\left\{\left\Vert
\int\limits_{t}^{t_2}B(Z)d\left({\bf W}_{t_1}-{\bf W}_{t_1}^M\right)
\right\Vert_H^4\right\}\right)^{1/2}dt_3\le
$$

\vspace{4mm}
$$
\le C_{13}\int\limits_{t}^{T}
\int\limits_{t}^{t_2}\left({\sf M}\left\{\biggl\Vert B(Z)
\biggr\Vert_{L_{HS}(U_0,H)}^4\right\}\right)^{1/2}
dt_1
\left({\sf M}\left\{\left\Vert
\int\limits_{t}^{t_2}B(Z)d\left({\bf W}_{t_1}-{\bf W}_{t_1}^M\right)
\right\Vert_H^4\right\}\right)^{1/2}dt_3\le
$$

\vspace{4mm}
\begin{equation}
\label{lll1}
\le C_{14}\int\limits_{t}^{T}
(t_2-t)
\left({\sf M}\left\{\left\Vert
\int\limits_{t}^{t_2}B(Z)d\left({\bf W}_{t_1}-{\bf W}_{t_1}^M\right)
\right\Vert_H^4\right\}\right)^{1/2}dt_3.
\end{equation}

\vspace{8mm}

Let us estimate the right-hand side of (\ref{lll1}).
Let $s>t.$ For fixed $M\in\mathbb{N}$ and 
for some $N>M$ ($N\in\mathbb{N}$) we have

\vspace{2mm}
$$
{\sf M}\left\{\left\Vert
\int\limits_{t}^{s}B(Z)d\left({\bf W}_{t_1}^N-{\bf W}_{t_1}^M\right)
\right\Vert_H^4\right\}=
$$

\vspace{4mm}
$$
={\sf M}\left\{
\Biggl\langle 
\sum\limits_{j\in J_N\backslash J_M}\sqrt{\lambda_j}
B(Z)e_j\left({\bf w}_{s}^{(j)}-{\bf w}_{t}^{(j)}\right),
\sum\limits_{j'\in J_N\backslash J_M}\sqrt{\lambda_{j'}}
B(Z)e_{j'}\left({\bf w}_{s}^{(j')}-{\bf w}_{t}^{(j')}\right)
\Biggr\rangle_H^2\right\}=
$$

\vspace{8mm}

$$
=
\sum\limits_{j, j', l, l'\in J_N\backslash J_M}\ 
\sqrt{\lambda_j\lambda_{j'}\lambda_l\lambda_{l'}}\ \
{\sf M}\left\{
\biggl\langle 
B(Z)e_j,
B(Z)e_{j'}
\biggr\rangle_H\biggl\langle 
B(Z)e_l,
B(Z)e_{l'}
\biggr\rangle_H\times\right.
$$

\vspace{4mm}

$$
\left.\times
{\sf M}\left\{\left({\bf w}_{s}^{(j)}-{\bf w}_{t}^{(j)}\right)
\left({\bf w}_{s}^{(j')}-{\bf w}_{t}^{(j')}\right)
\left({\bf w}_{s}^{(l)}-{\bf w}_{t}^{(l)}\right)
\left({\bf w}_{s}^{(l')}-{\bf w}_{t}^{(l')}\right)
\biggl.\biggr|{\bf F}_t\right\}\right\}
=
$$

\vspace{8mm}

$$
=3(s-t)^2\sum\limits_{j\in J_N\backslash J_M}
\lambda_j^{2}
{\sf M}\left\{\biggl\Vert B(Z)e_j\biggr\Vert_H^4\right\}+
$$

\vspace{4mm}
$$
+(s-t)^2
\sum\limits_{j, j'\in J_N\backslash J_M (j\ne j')}
\lambda_j\lambda_{j'}\Biggl(
{\sf M}\left\{\biggl\Vert B(Z)e_j\biggr\Vert_H^2
\biggl\Vert B(Z)e_{j'}\biggr\Vert_H^2\right\}\Biggr.
\Biggl.+
2\biggl\langle 
B(Z)e_j,
B(Z)e_{j'}
\biggr\rangle_H^2\Biggr)\le
$$

\vspace{8mm}

$$
\hspace{-60mm}
\le 3(s-t)^2\left(\sum\limits_{j\in J_N\backslash J_M}
\lambda_j^{2}
{\sf M}\left\{\biggl\Vert B(Z)e_j\biggr\Vert_H^4\right\}+\right.
$$

\vspace{1mm}

$$
~~~~~~~~~~~~~~~~~~~~~~~~~~~~~~~~~~~~~~~~~+\left.
\sum\limits_{j, j'\in J_N\backslash J_M (j\ne j')}
\lambda_j\lambda_{j'}
{\sf M}\left\{\biggl\Vert B(Z)e_j\biggr\Vert_H^2
\biggl\Vert B(Z)e_{j'}\biggr\Vert_H^2\right\}\right)=
$$

\vspace{6mm}

$$
=
3(s-t)^2
{\sf M}\left\{\left(
\sum\limits_{j\in J_N\backslash J_M}\lambda_j
\biggl\Vert B(Z)e_j\biggr\Vert_H^2\right)^2\right\}\le
$$

\vspace{4mm}

$$
\le
3(s-t)^2\left(\sup\limits_{i\in J_N\backslash J_M}\lambda_i\right)^{4\alpha}
{\sf M}\left\{\left(
\sum\limits_{j\in J_N\backslash J_M}\lambda_j
\biggl\Vert B(Z)Q^{-\alpha}e_j\biggr\Vert_H^2\right)^2\right\}\le
$$

\vspace{4mm}
\begin{equation}
\label{hhhh1}
\le
C_{15}(s-t)^2
\left(\sup\limits_{i\in J_N\backslash J_M}\lambda_i\right)^{4\alpha}
{\sf M}\left\{
\biggl\Vert B(Z)Q^{-\alpha}\biggr\Vert_{L_{HS}(U_0,H)}^4\right\}.
\end{equation}

\vspace{10mm}

Performing the passage to the 
limit $\lim\limits_{N\to\infty}$ in (\ref{hhhh1})
and 
using (\ref{lll1}), we have

\vspace{2mm}
\begin{equation}
\label{aaa1}
G^{1,M}_{T,t}\le C_{16} 
\left(\sup\limits_{i\in J\backslash J_M}\lambda_i\right)^{2\alpha}
(T-t)^3.
\end{equation}

\vspace{6mm}

Absolutely analogously we obtain

\vspace{2mm}
\begin{equation}
\label{aaa2}
G^{2,M}_{T,t}\le C_{17} 
\left(\sup\limits_{i\in J\backslash J_M}\lambda_i\right)^{2\alpha}
(T-t)^3.
\end{equation}

\vspace{6mm}

Let us estimate $G^{3,M}_{T,t}.$ We have

\vspace{2mm}
$$
G^{3,M}_{T,t}\le
\left(\sup\limits_{i\in J\backslash J_M}\lambda_i\right)^{2\alpha}\times
$$

\vspace{4mm}

$$
\times
\int\limits_{t}^{T}
{\sf M}\left\{\left\Vert
B''(Z) \left(
\int\limits_{t}^{t_2}B(Z)d{\bf W}_{t_1}^M,
\int\limits_{t}^{t_2}B(Z)d{\bf W}_{t_1}^M
\right)Q^{-\alpha} 
\right\Vert_{L_{HS}(U_0,H)}^2\right\}dt_2\le
$$

\vspace{4mm}

$$
\le \left(\sup\limits_{i\in J\backslash J_M}\lambda_i\right)^{2\alpha}
\sum\limits_{i\in J}
\sum\limits_{j,l\in J_M}\lambda_i \lambda_j
\lambda_l
\int\limits_t^T(t_2-t)^2 \times
$$

\vspace{4mm}
$$
\times \Biggl( {\sf M} \Biggl\{
\Biggl\Vert B''(Z)(B(Z)e_j,B(Z)e_l)Q^{-\alpha}e_i
\Biggr\Vert_H^2\Biggr\} +\Biggr.
$$

\vspace{4mm}

$$
\hspace{-40mm}
+ {\sf M} \left\{\Biggl\Vert B''(Z)(B(Z)e_j,B(Z)e_j)Q^{-\alpha}e_i
\Biggr\Vert_H\times\right.
$$

\vspace{4mm}

$$
\left.~~~~~~~~~~~~~~~~~~~~~~~~~~~~~~~~~~~~~~~~\times
\Biggl\Vert B''(Z)(B(Z)e_l,B(Z)e_l)Q^{-\alpha}e_i
\Biggr\Vert_H\right\}
+
$$

\vspace{4mm}

$$
\hspace{-40mm}
+
{\sf M} \left\{\Biggl\Vert B''(Z)(B(Z)e_j,B(Z)e_l)Q^{-\alpha}e_i
\Biggr\Vert_H\times\right.
$$

\vspace{4mm}
$$
\left.\left.~~~~~~~~~~~~~~~~~~~~~~~~~~~~~~~~~~~~~~~~\times
\Biggl\Vert B''(Z)(B(Z)e_l,B(Z)e_j)Q^{-\alpha}e_i
\Biggr\Vert_H
\right\}\right)dt_2\le
$$

\vspace{4mm}
\begin{equation}
\label{kkk1}
\le C_{18}
\left(\sup\limits_{i\in J\backslash J_M}\lambda_i\right)^{2\alpha}
(T-t)^3.
\end{equation}

\vspace{8mm}

Combaining (\ref{sss1}), (\ref{ppp0x}), and (\ref{aaa1})--(\ref{kkk1}),
we get (\ref{uuu1}). Theorem 4 is proved.

Let us consider the convergence analysis
for the stochastic integrals
(\ref{ret2})--(\ref{ret4}) (convergence of the 
stochastic integral (\ref{ret1}) follows from 
(\ref{24u}) (see Theorem 1 in \cite{8} or \cite{2})). 

Using the It\^{o} formula, we obtain w. p. 1 \cite{3}

\vspace{3mm}
$$
J_2[B(Z)]_{T,t}=\int\limits_{t}^{T}\left(\frac{T}{2}-s+\frac{t}{2}\right)
AB(Z)d{\bf W}_{s},
$$

\vspace{2mm}
$$
J_3[B(Z),F(Z)]_{T,t}
=\int\limits_{t}^{T}(s-t)B'(Z)\biggl(AZ+F(Z)\biggr)d{\bf W}_{s}.
$$

\vspace{6mm}

Suppose that

\vspace{2mm}
$$
{\sf M}\left\{\biggl\Vert B'(Z)\biggl(AZ+F(Z)\biggr)Q^{-\alpha}
\biggr\Vert_{L_{HS(U_0,H)}}^2\right\}<\infty,
$$

\vspace{4mm}
$$
{\sf M}\left\{\biggl\Vert AB(Z)Q^{-\alpha}\biggr\Vert_{L_{HS(U_0,H)}}^2\right\}
<\infty
$$

\vspace{6mm}
\noindent
for some $\alpha>0.$

Then by analogy with (\ref{24u}) we get

\vspace{2mm}
$$
{\sf M}\left\{\Biggl\Vert
J_2[B(Z)]_{T,t}-
J_2[B(Z)]_{T,t}^{M}\Biggr\Vert_H^2\right\}\le
$$

\vspace{2mm}
$$
\le C_{19} (T-t)^3
\left(\sup\limits_{i\in J\backslash J_M}\lambda_i\right)^{2\alpha},
$$

\vspace{7mm}
$$
{\sf M}\left\{\Biggl\Vert
J_3[B(Z),F(Z)]_{T,t}-
J_3[B(Z),F(Z)]_{T,t}^{M}\Biggr\Vert_H^2\right\}\le
$$

\vspace{2mm}
$$
\le C_{20} (T-t)^3
\left(\sup\limits_{i\in J\backslash J_M}\lambda_i\right)^{2\alpha},
$$

\vspace{9mm}
\noindent
where
$J_2[B(Z)]_{T,t}^{M},$
$J_3[B(Z),F(Z)]_{T,t}^{M}$
are defined by 
(\ref{sd11}), (\ref{sd12}).

Moreover, in conditions of Theorem 4 we obtain

\vspace{3mm}
$$
{\sf M}\left\{\Biggl\Vert
J_4[B(Z),F(Z)]_{T,t}-
J_4[B(Z),F(Z)]_{T,t}^{M}\Biggr\Vert_H^2\right\}=
$$

\vspace{5mm}
$$
={\sf M}\left\{\left\Vert
\int\limits_t^T F'(Z)\left(\int\limits_t^{t_2}
B(Z)
d\left({\bf W}_{t_1}-{\bf W}_{t_1}^M\right)\right)
dt_2\right\Vert_H^2\right\}\le
$$

\vspace{5mm}
$$
\le 
(T-t)
\int\limits_t^T
{\sf M}\left\{\left\Vert
F'(Z)\left(\int\limits_t^{t_2}
B(Z)
d\left({\bf W}_{t_1}-{\bf W}_{t_1}^M\right)\right)
\right\Vert_H^2\right\}dt_2\le
$$

\vspace{5mm}
$$
\le 
C_{21}(T-t)
\int\limits_t^T
{\sf M}\left\{\left\Vert
\int\limits_t^{t_2}
B(Z)
d\left({\bf W}_{t_1}-{\bf W}_{t_1}^M\right)
\right\Vert_H^2\right\}dt_2\le
$$

\vspace{5mm}
$$
\le C_{21}(T-t)
\left(\sup\limits_{i\in J\backslash J_M}\lambda_i\right)^{2\alpha}
\int\limits_t^T
\int\limits_t^{t_2}
{\sf M}\left\{\biggl\Vert
B(Z)Q^{-\alpha}
\biggr\Vert_{L_{HS}(U_0,H)}^2\right\} dt_1dt_2\le
$$

\vspace{4mm}
$$
\le 
C_{22}(T-t)^3
\left(\sup\limits_{i\in J\backslash J_M}\lambda_i\right)^{2\alpha},
$$

\vspace{9mm}
\noindent
where 
$J_4[B(Z),F(Z)]_{T,t}^{M}$ is defined by 
(\ref{sd13}).

\end{document}